\documentclass[a4paper,12pt]{amsart}
\usepackage[english]{babel}
\usepackage[T1]{fontenc}

\usepackage{enumerate}
\usepackage{comment}
\usepackage{multirow}

\usepackage[a4paper,margin=1.1in]{geometry}

\usepackage{amsfonts,amsmath,amsthm,amssymb,amsxtra,calligra,mathrsfs}
\usepackage{graphicx}
\usepackage{tikz}
\usepackage{tikz-cd} 
\usepackage{xcolor}
\usepackage{upgreek}
\usepackage{graphicx}
\usepackage{mathtools}
\usepackage{extarrows}
\usepackage{faktor}
\usepackage{tensor}
\usepackage{enumitem}
\usepackage{hyperref}
\hypersetup{
    colorlinks=true, 
    citecolor=blue,
    linkcolor=magenta,
	bookmarksnumbered=true,
}

\newcommand{\bZ}{\mathbb{Z}}
\newcommand{\bP}{\mathbb{P}}
\newcommand{\bC}{\mathbb{C}}

\newcommand{\bQ}{{\mathbb Q}}

\newcommand{\bH}{{\mathbb H}}

\newcommand{\cO}{{\mathcal O}}

\newcommand{\cZ}{{\mathcal Z}}

\newcommand{\cI}{{\mathcal I}}
\newcommand{\cH}{{\mathcal H}}
\newcommand{\cG}{{\mathcal G}}
\newcommand{\cA}{{\mathcal A}}

\newcommand{\sheafHom}{\mathscr{H}\text{\kern -3pt {\calligra\large om}}\,}

\def\D{\Delta}

\DeclareMathOperator{\IH}{\mathrm{IH}}

\DeclareMathOperator{\Spec}{Spec}

\DeclareMathOperator{\Supp}{Supp}
\DeclareMathOperator{\rk}{rk}
\DeclareMathOperator{\clg}{Cl}

\DeclareMathOperator{\Ext}{Ext}
\DeclareMathOperator{\cext}{\mathcal{E}xt}

\DeclareMathOperator{\Pic}{\mathrm{Pic}}
\DeclareMathOperator{\NS}{\mathrm{NS}}

\DeclareMathOperator{\CH}{\mathrm{CH}}

\DeclareMathOperator{\ch}{\mathrm{ch}}

\DeclareMathOperator{\codim}{\mathrm{codim}}

\DeclareMathOperator{\pic}{Pic}

\newcommand\im{\text{\rm Im}}
\newcommand\rank{{\text{\rm rank}}}
\newcommand\lra{\longrightarrow}

\newcommand{\cF}{{\mathcal F}}

\newcommand\cE{{\mathcal{E}}}

\newcommand\cL{\mathcal{L}}

\newcommand\cC{{\mathcal{C}}}

\newcommand\cK{{\mathcal{K}}}

\newtheorem{defn-pro}{Definition-Proposition}
\newtheorem{defn-thm}{Definition-Theorem}
\newtheorem{thm}{Theorem}[section]
\newtheorem{assu}[thm]{Assumption}
\newtheorem{lem}[thm]{Lemma}
\newtheorem{cor}[thm]{Corollary}
\newtheorem{pro}[thm]{Proposition}

\newtheorem{defn}[thm]{Definition}

\newtheorem{question}[thm]{Question}
\newtheorem{exmp}[thm]{Example}

\newtheorem{rem}[thm]{Remark}

\theoremstyle{remark}

\title[Picard numbers of Moduli of One-dimensional Sheaves on Surfaces]{On the Picard numbers of Moduli spaces of One-dimensional Sheaves on Surfaces}

\author{Fei Si}
\address{The School of Mathematics and Statistics, Xi’an Jiaotong University, 28 West Xianning Road, Xi’an, Shaanxi, P.R.China 710049}
\email{sifei@xjtu.edu.cn}

\author{Feinuo Zhang}
\address{Shanghai Center for Mathematical Sciences, Fudan University, Jiangwan Campus, Shanghai, 200438, China}
\email{fnzhang17@fudan.edu.cn}
\date{}

\usepackage{microtype}
\begin{document}

\begin{abstract}
  Motivated by asymptotic phenomena of moduli spaces of higher rank stable sheaves on algebraic surfaces, we study the Picard number of the moduli space of one-dimensional stable sheaves supported in a sufficiently positive divisor class on a surface. We give an asymptotic lower bound of the Picard number in general. In some special cases, we show that this lower bound is attained based on the geometry of moduli spaces of stable pairs and relative Hilbert schemes of points. Additionally, we discuss several related questions and provide examples where the asymptotic irreducibility of the moduli space fails, highlighting a notable distinction from the higher rank case. 
\end{abstract}
\maketitle

\tableofcontents

\section{Introduction}

Let $S$ be a smooth, complex projective surface with an ample divisor $H$ (such a pair $(S,H)$ is called a polarized surface). 
The moduli spaces of sheaves on $S$ have played a crucial role in various areas of mathematics, and their geometry has long fascinated algebraic geometers. We denote by $M(r,L,c_2)$ the (coarse) moduli space parametrizing polystable sheaves (with respect to $H$) on $S$ of prescribed rank $r\in\bZ_{\geq0}$, with fixed determinant $L$ and second Chern class $c_2\in H^4(S,\bZ)\cong\bZ$.

When $r=1$, $M(1,L,c_2)$ is isomorphic to $S^{[c_2]}$, the Hilbert scheme of $0$-dimensional closed subschemes in $S$ of length $c_2$. Fogarty showed that $S^{[n]}$ ($n\in\bZ_{\geq0}$) is a smooth projective variety of dimension $2n$ \cite{Fog68} and computed its Picard group \cite{Fog73}. The Betti numbers of $S^{[n]}$ were completely determined by Göttsche \cite{Got90}. 
When $r\geq 2$, $M(r,L,c_2)$ can be reducible or non-reduced. But starting from Donaldson's generic smoothness result \cite{Don90}, the study of $M(r,L,c_2)$ has revealed that $M(r,L,c_2)$ has better geometric properties as we fix $r\geq2$, $L$, and let $c_2$ increase. For example, when $c_2$ is large enough, Gieseker--Li \cite{GL94} \cite{GiesekerLi96} and O'Grady \cite{O'Grady96} showed that $M(r,L,c_2)$ is irreducible, normal, and generically smooth of expected dimension.
Some asymptotic phenomena also occur in the Picard group and Betti numbers of $M(r,L,c_2)$ ($r\geq2$) \cite{Qin92}\cite{Li96}\cite{Li97}\cite{CW22}.

This raises a natural question:
\begin{question}
\label{quest}
Are there asymptotic phenomena similar to those above in the geometry of $M(r,L,c_2)$ when $r=0$?
\end{question}

Sheaves on $S$ of rank zero are supported on points or curves. In this paper, we focus on those pure sheaves supported on curves, which are called one-dimensional sheaves. In contrast to the positive rank case, there are few general results for the moduli space of one-dimensional sheaves. An essential difficulty is that the supports of one-dimensional sheaves can be singular, even non-reduced. To distinguish this case, we use the following notation.
Given a nonzero effective divisor $\beta$ on $S$ and an integer $\chi$, denote by $M_{\beta,\chi}$ the moduli space parametrizing polystable one-dimensional sheaves $F$ on $S$, with determinant $\det(F)\cong\cO_S(\beta)$ and Euler characteristic $\chi(F)=\chi$. 
We will make the following assumption on $\beta$, unless otherwise stated.
\begin{assu}
\label{assu:beta}
The following conditions hold:
    \begin{enumerate}
    \item $\beta$ is base-point-free.
    \item A general curve in $|\beta|$ is smooth and connected.
\end{enumerate}
\end{assu}

A conjecture \cite[Conjecture 1.3]{SZ24} was proposed on the stabilization of (intersection cohomology) Betti numbers of $M_{\beta,\chi}$ when $\chi$ is fixed and $\beta$ is sufficiently positive, as an analogue of \cite[Conjecture 1.1]{CW22}. When $S$ is a del Pezzo surface, this conjecture was resolved in \cite{PSSZ24} by proving an asymptotic formula of the refined BPS invariants which also refine the Betti numbers. 

The primary motivation of this paper is to explore the stabilization of the Picard number $\rho(M_{\beta,\chi})$ of $M_{\beta,\chi}$ when $\beta$ is sufficiently positive. Here, the Picard number $\rho(X)$ of a proper scheme $X$ is defined as the dimension of the $\bQ$-vector space $\mathrm{Num}_\bQ(X)$ of the numerically equivalence classes of line bundles on $X$ with rational coefficients, in the sense of \cite{Kle66}. Fix a basis $[L_i]$ $(i=1,\cdots, \rho(S))$ for
$\mathrm{Num}_{\bQ}(S)$, 
where each $L_i$ is a smooth, connected effective divisor. We specify a positivity condition for $\beta$ as follows.

\begin{defn}
\label{P}
    We say a divisor $\beta$ on $S$ satisfies $(\mathrm{P})$ if
    \begin{enumerate}
    \item $\beta$ satisfies Assumption \ref{assu:beta}.
    \item The arithmetic genus $p_a(\beta)>0$.
        \item For each $1\leq i\leq\rho(S)$, there exists a smooth and connected curve $L_i'\in|\beta-L_i|$ such that $L_i$ meets $L_i'$ transversally with intersection number $L_i\cdot L_i'\geq|\chi|+1$.
    \end{enumerate}
\end{defn}

Note that for an ample divisor $\beta$, the multiple $n\beta$ satisfies $(\mathrm{P})$ for all $n\gg0$. 

Our first main result is the following. 
We obtain an asymptotic lower bound of the Picard number of $M_{\beta,\chi}$.

\begin{thm}[Theorem \ref{picardnum}]
\label{mainpic}
Let $(S,H)$ be any polarized surface, and let $\beta$ be a nonzero effective divisor satisfying $(\mathrm{P})$ as defined in Definition \ref{P} such that $\beta\cdot H$ and $\chi\in\bZ$ are coprime.
Then we have 
\begin{equation}
\label{rho}
    \rho(M_{\beta,\chi})\geq\rho(S)+1.
\end{equation}
\end{thm}

The above result generalizes \cite[Theorem 9.2]{CW22} to the rank zero case.
We prove it by using the determinant line bundles and constructing testing curves in $M_{\beta,\chi}$. In our case, it is trickier to construct testing curves than in the positive rank case since Hecke modifications at points may change the stability of a one-dimensional sheaf. The condition $(\mathrm{P})$ is needed to ensure the existence of certain nodal curves in $|\beta|$. 

\begin{question}
\label{ques:pic}
In the situation of Theorem \ref{mainpic}, if $h^1(\cO_S):=\dim_\bC H^1(\cO_S)=0$, does $\rho(M_{\beta,\chi})$ stabilize to $\rho(S)+1$ when $\beta$ is sufficiently positive?
\end{question}

This question is answered affirmatively in the case of K3 surfaces using the hyperkähler geometry of $M_{\beta,\chi}$ and \cite[Corollary 6.3]{Fog73}, and in the case of del Pezzo surfaces by \cite{PSSZ24}. In fact, Theorem \ref{mainpic} combined with the tautological generation result \cite{markman07} recovers these two cases. Another evidence for the stabilization of $\rho(M_{\beta,\chi})$ is given by \cite{Sac19} in the case of Enriques surfaces, where tautological generation fails. In this paper, we also provide evidence for a positive answer to Question \ref{ques:pic} beyond the above cases under some assumptions.
\begin{thm}
\label{mt:piceq}
Let $(S,H)$ be any polarized surface with $h^1(\cO_S)=0$, and let $\beta$ be a nonzero effective divisor satisfying Assumption \ref{assu:beta} such that $\beta\cdot H$ and $\chi>1$ are coprime.
Suppose $M_{\beta,\chi}$ is normal and irreducible, with the property that the locus of sheaves with non-integral support has codimension at least $2$. In the situations (1)-(3) of Theorem \ref{thm4}, if $\beta$ also satisfies $(\mathrm{P})$ as defined in Definition \ref{P}, we have $\rho(M_{\beta,\chi})=\rho(S)+1$.
\end{thm}

We emphasize that the geometry of moduli spaces of one-dimensional sheaves is not a simple generalization of the positive rank case, since the techniques for torsion free sheaves are usually not applicable to torsion sheaves, which makes the above questions more interesting. Indeed, we discover the following phenomenon for those irregular surfaces whose canonical divisor $K_S$ is quite positive (e.g., $K_S$ is very ample).
\begin{thm}[Theorem \ref{counter}]
Let $(S,H)$, together with an integer $n\geq 2$, be a polarized irregular surface satisfying Assumption \ref{assu:irrd}. If $\beta_n=nK_S$, then the moduli space $M_{\beta_n,1}$ is reducible.
\end{thm}

Our examples show that when $r=0$, the asymptotic irreducibility of moduli spaces of one-dimensional sheaves can fail, which is in nature different from the higher rank case proved by Gieseker--Li and O'Grady.
However, this phenomenon does not happen when $S$ is a del Pezzo surface (see, e.g., \cite[Theorem 2.3]{MS23}\cite[Theorem 1.5]{Yuan23}), suggesting an intricate pattern in general.

The rest of this paper is organized as follows. In Section \ref{sec:pre}, we collect some notation, definitions, and facts which will be useful later. In Section \ref{sec:pic}, we construct three types of testing curves and prove the asymptotic lower bound for the Picard number. The approach to an upper bound motivated by the wall-crossing of moduli spaces of pairs is discussed in Section \ref{sec:wc}. In the end, we present in Section \ref{sec:other} examples in which the moduli spaces are not irreducible.

\subsection*{Conventions}
Throughout this paper, we work over the field $\bC$ of complex numbers. 
All schemes are assumed to be of finite type over $\bC$. Points of a scheme are assumed to be closed. By a variety, we mean an integral, separated scheme. By a surface, we mean a smooth, projective $2$-dimensional variety. A curve is a projective scheme of pure dimension $1$. All sheaves are assumed to be coherent. For a divisor $D$ on a smooth projective variety $X$, we denote by $|D|$ the complete linear system $\bP(H^0(\cO_X(D))^\vee)$.

\subsection*{Acknowledgments}
The second author is grateful to her advisor, Jun Li for his help in preparing for this work, and to Claire Voisin for useful discussions. Both authors thank Zhiyuan Li, Weite Pi, and Junliang Shen for helpful comments.
The first author is supported by NSFC grant 12201011. The second author is supported by NKRD Program of China (No. 2020YFA0713200) and LMNS.

\section{Preliminaries}
\label{sec:pre}
In this section, we review some properties concerning moduli spaces of our interest.

\subsection{Moduli of one-dimensional sheaves and determinant line bundles}
Let $(S,H)$ be a polarized surface. For a nonzero sheaf $E$ on $S$, the Hilbert polynomial of $E$ with respect to $H$ is
$$\chi(E\otimes \cO_S(mH))=\sum_{i=0}^{d_E}a_i(H,E)\frac{m^i}{i!}\quad (m\in\bZ_{\geq0}),$$
whose degree $d_E\in\bZ_{\geq0}$ is exactly the dimension of the support of $E$, called the dimension of $E$.
The sheaf $E$ is called one-dimensional if every nonzero subsheaf of $E$ has dimension $1$. In this case, the slope $\mu(E)$ is defined by
$$\mu(E):=\frac{\chi(E)}{c_1(E)\cdot H}.$$

\begin{defn}[\cite{Simpson94}]
\label{def:stab}
   A one-dimensional sheaf $F$ on $S$ is called semistable (with respect to $H$) if for any nonzero proper subsheaf $G\subset F$, 
   \begin{equation}
   \label{ssineq}
       \mu(G) \le  \mu(F).
   \end{equation}
   The above sheaf $F$ is called stable if it is semistable and the inequality (\ref{ssineq}) is strict. A semistable sheaf is called polystable if it is the direct sum of stable sheaves.
\end{defn}

For any one-dimensional sheaf $F$, any surjection $f:E_0\to F$ from a locally free sheaf $E_0$ can be completed to a short exact sequence
$$0\to E_1\stackrel{g}\to E_0\stackrel{f}\to F\to0,$$
where $E_1$ is locally free of rank $r(E_1)=r(E_0)$. The Fitting support $\Supp(F)$ of $F$ is defined (as a closed subscheme of $S$) by the vanishing locus of the determinant of $g$, which is independent of the locally free resolution $E_{\bullet}$ of $F$  \cite[\S 20.2, p. 493]{Eis95}.
Given a nonzero effective divisor $\beta$ and $\chi\in\bZ$, the moduli space $M_{\beta,\chi}$ admits a Hilbert--Chow morphism 
$$h_{\beta,\chi}:M_{\beta,\chi} \rightarrow |\beta|$$ 
which sends a polystable sheaf $F\in M_{\beta,\chi}$ to its Fitting support $\Supp(F)$. 



Let $U\subset|\beta|$ be the open subset parametrizing integral curves.
Then there is another modular interpretation of $h_{\beta,\chi}^{-1}(U)$ as the relative compactified Jacobian of degree $d$ associated to the universal family  $\pi:\cC_U\rightarrow U$ 
of integral curves in $|\beta|$, where $\cC_U \subset U \times S$ is a closed subscheme  and  
\begin{equation*}
    d=\chi+\frac{\beta(\beta+K_S)}{2}. 
\end{equation*}

\begin{pro}
\label{prop:smsupp}
The open subscheme $h_{\beta,\chi}^{-1}(U)$ of $M_{\beta,\chi}$ is irreducible and contains a smooth open subscheme $h_{\beta,\chi}^{-1}(U_0)$ of dimension
\begin{equation}
    \label{dimintsm}
    \dim h_{\beta,\chi}^{-1}(U_0)=\dim h_{\beta,\chi}^{-1}(U)=h^0(\cO_S(\beta))+\frac{K_S\cdot\beta+\beta^2}{2},
\end{equation}
where $U_0\subset|\beta|$ is the locus of smooth curves. If $h^1(\cO_S(\beta))=h^2(\cO_S(\beta))=0$, then 
     \begin{equation}
         \label{dimsmsupp}
          \dim h_{\beta,\chi}^{-1}(U_0)=\dim h_{\beta,\chi}^{-1}(U)=\beta^2+\chi(\cO_S).
     \end{equation}   
\end{pro}
\begin{proof}
The fibers of $h_{\beta,\chi}$ over $U$ are integral of dimension $p_a(\beta)$ by \cite[Theorem (9)]{AIK76}, so $h_{\beta,\chi}^{-1}(U)$ is irreducible of dimension $\dim|\beta|+p_a(\beta)$.
    Since semistable sheaves on smooth curves with Euler characteristic $\chi$ are line bundles of degree $d$, the open subscheme $h_{\beta,\chi}^{-1}(U_0)$ is the relative Picard scheme $\Pic^d(\cC_{U_0}/U_0)$, where $\cC_{U_0}:=\pi^{-1}(U_0)$. By \cite[Corollary 5.14, Proposition 5.19]{Kleiman05}, $\Pic^d(\cC_{U_0}/U_0)$ is smooth over $U_0$ of relative dimension $p_a(\beta)$. Thus $h_{\beta,\chi}^{-1}(U_0)$ is smooth and (\ref{dimintsm}) follows by the adjunction formula. 
    
    If $h^1(\cO_S(\beta))=h^2(\cO_S(\beta))=0$, then by 
    the Riemann--Roch formula 
    \begin{equation*}
        \begin{aligned}
    \dim h_{\beta,\chi}^{-1}(U_0)&= \dim h_{\beta,\chi}^{-1}(U)=\chi(\cO_S(\beta))+\frac{K_S\cdot\beta+\beta^2}{2}\\
    &=\left(\frac{\beta(\beta-K_S)}{2}+\chi(\cO_S)\right)+\frac{K_S\cdot\beta+\beta^2}{2}\\
    &=\beta^2+\chi(\cO_S),
        \end{aligned}
    \end{equation*}
which completes the proof.
\end{proof}


Now we assume that $\beta\cdot H$ and $\chi$ are coprime. In this case, all sheaves in $M_{\beta,\chi}$ are stable. 
By \cite[Corollary 4.6.6]{HL10}, there is a universal sheaf $\cE$ on $M_{\beta,\chi}\times S$.
Recall that the Grothendieck group $K^0(X)$ of locally free sheaves on a scheme $X$ is a commutative ring with $1=[\cO_X]$ and the multiplication given by 
$$[F_1]\cdot[F_2]:=[F_1\otimes F_2]$$ for locally free sheaves $F_1$ and $F_2$ on $X$.
There is a group homomorphism from $K^0(S)$ to the Picard group of $M_{\beta,\chi}$ (\cite[Definition 8.1.1]{HL10})
\begin{equation*}
   \lambda_{\cE}:\  K^0(S) \rightarrow \pic(M_{\beta,\chi}),\quad  \lambda_\cE(\mathfrak{a}):=\det(p_{!}([\cE] \cdot q^\ast \mathfrak{a} ))\quad(\mathfrak{a}\in K^0(S)),
\end{equation*}
where $p:M_{\beta,\chi}\times S\to M_{\beta,\chi}$ and $q:M_{\beta,\chi}\times S\to S$ are the projections. Note that the class $[\cE]\in K^0(M_{\beta,\chi}\times S)$ and the homomorphism $p_!:K^0(M_{\beta,\chi}\times S)\to K^0(M_{\beta,\chi})$ are well-defined by \cite[Corollary 2.1.11]{HL10}.

Let $K_{\beta,\chi}^0\subset K^0(S)$ be the orthogonal complement of the class of a sheaf in $M_{\beta,\chi}$ with respect to the bilinear form $(\mathfrak{a},\mathfrak{b})\mapsto \chi(\mathfrak{a}\cdot \mathfrak{b})$ on $K^0(S)\times K^0(S)$.
Since two universal sheaves differ by a line bundle pulled back from $M_{\beta,\chi}$, the homomorphism $\lambda_\cE$ is independent of the choice of $\cE$ when restricted to $K_{\beta,\chi}^0$ by \cite[Lemma 8.1.2 {\it iv)}]{HL10}. We denote the restriction by
\begin{equation}
    \label{lambda}
    \lambda:K_{\beta,\chi}^0\to \pic(M_{\beta,\chi}).
\end{equation}

For a sheaf $\cF$ on a scheme $X$ and a sheaf $\cG$ on a scheme $Y$, we write $\cF\boxtimes \cG$ for the sheaf $p_X^*\cF\otimes p_Y^*\cG$ on $X\times Y$, where $p_X:X\times Y\to X$ and $p_Y:X\times Y\to Y$ are the projections.

\begin{pro}\label{firstchern}
The determinant of $[\cE]\in K^0(M_{\beta,\chi}\times S)$ is
\[\det(\cE)=h_S^\ast (\cO_{|\beta|}(1) \boxtimes \cO_S(\beta)),   \] 
 where $h_S:=h_{\beta,\chi}\times 1_S: M_{\beta,\chi} \times S \rightarrow |\beta| \times S$.
\end{pro}
\begin{proof}
  The argument is similar to that of \cite[Lemma 2.1]{PS23}. Clearly, $\cE$ is a torsion sheaf (on $M_{\beta,\chi} \times S$) supported on $h^{-1}(\cC)$, the pullback of the universal curve $\cC\subset|\beta| \times S$ along $h$. Note that $\cC$ is flat over $|\beta|$ and therefore $h^{-1}(\cC)$ is a Cartier divisor on $M_{\beta,\chi}\times S$ corresponding to the pullback of the line bundle
  \begin{equation*}
      \cO_{|\beta| \times S}(\cC) \cong \cO_{|\beta|}(1) \boxtimes \cO_S(\beta),
  \end{equation*}
  from which the result follows.
\end{proof}

\begin{cor}\label{cor:lambda0}
    Let $x\in S$ be a point and denote by $[\cO_x]\in K^0(S)$ the class of the structure sheaf of $x$. Then $[\cO_x]\in K_{\beta,\chi}^0$ and
    $$\lambda([\cO_x])=h_{\beta,\chi}^*\cO_{|\beta|}(1).$$
\end{cor}
\begin{proof}
The assertion that $[\cO_x]\in K_{\beta,\chi}^0$ follows immediately from a locally free resolution of $F\in M_{\beta,\chi}$.
    By \cite[Example 8.1.3 {\it i)}]{HL10},
    $$\lambda([\cO_x])=p_*(\det(\cE)|_{M_{\beta,\chi}\times\{x\}}).$$
    It follows from Proposition \ref{firstchern} that
    $$p_*(\det(\cE|_{M_{\beta,\chi}\times\{x\}}))=p_*(p^*h_{\beta,\chi}^*\cO_{|\beta|}(1))=h_{\beta,\chi}^*\cO_{|\beta|}(1),$$
    which completes the proof.
\end{proof}

\subsection{Relative Hilbert schemes of points}
\label{subsec:relHil}
For $n\in\bZ_{\geq0}$, denote by $\pi^{[n]}:\cC^{[n]}\to |\beta|$ the relative Hilbert scheme of $n$ points on the universal family $\cC$ of curves in $|\beta|$. For $C\in|\beta|$, the fiber of $\pi^{[n]}$ over $C$ is the Hilbert scheme $C^{[n]}$ parametrizing length $n$, zero-dimensional closed subschemes of $C$. Note that $\pi^{[0]}:|\beta|\to|\beta|$ is the identity and
$\pi:=\pi^{[1]}:\cC\to|\beta|$ is the natural projection.

The following notion of being $k$-very ample is a generalization of being base-point-free ($0$-very ample) and being very ample ($1$-very ample).

\begin{defn}
    A nonzero effective divisor $\beta$ on $S$ is called $k$-very ample $(k\in\bZ_{\geq0})$ if given any zero-dimensional closed subscheme $Z\subset S$ of length $k+1$, the restriction map 
    $$r_Z:H^0(\cO_S(\beta))\to H^0(\cO_S(\beta)|_Z)$$ is surjective.
\end{defn}



Since all the curves in $|\beta|$ lie in $S$, there is a natural morphism 
\begin{equation}
    \label{eq:defsigma}
    \sigma_k:\cC^{[k+1]}\to S^{[k+1]}
\end{equation}
that sends a closed subscheme $Z\subset C$ $(C\in|\beta|)$ of length $k+1$ to $Z\subset S$. The morphism $\sigma_k$ allows us to regard $\cC^{[k+1]}$ as a scheme over $S^{[k+1]}$. The following property of $\sigma_k$ is well-known (e.g., \cite[Proposition 3.16]{CGKT20}). We give a quick proof for the reader's convenience.

\begin{pro}
\label{kveryamp}
    If $\beta$ is $k$-very ample, then $\cC^{[k+1]}$ is a projective bundle over $S^{[k+1]}$.
\end{pro}
\begin{proof}
Let $\cI$ be the universal ideal sheaf on $S^{[k+1]}\times S$. Denote by $p_k:S^{[k+1]}\times S\to S^{[k+1]}$ and $p_S:S^{[k+1]}\times S\to S$ the projections. For every zero-dimensional closed subscheme $Z\subset S$ of length $k+1$ defined by the ideal sheaf $\cI_Z$, we have a short exact sequence
$$0\lra H^0(\cO_S(\beta)\otimes\cI_Z)\lra H^0(\cO_S(\beta))\stackrel{r_Z}\lra H^0(\cO_S(\beta)|_Z)\lra 0$$
since $\beta$ is $k$-very ample, which implies that
$h^0(\cO_S(\beta)\otimes\cI_Z)=h^0(\cO_S(\beta))-(k+1)$ is independent of the choice of $Z$.
By Grauert's theorem \cite[Ch. III, Corollary 12.9]{Hartshorne}, 
$$\cG:=p_{k\,*}(\cI\otimes p_S^*\,\cO_S(\beta))$$ is a vector bundle on $S^{[k+1]}$ of rank $h^0(\cO_S(\beta))-(k+1)$. Then it is easy to show that $\cC^{[k+1]}\cong  \bP(\cG^\vee)$.
\end{proof}

We conclude this section with the following equidimensional result, which will be useful in Section \ref{sec:wc}.
\begin{pro}
\label{eqdim:cC}
    The morphism $\pi^{[n]}:\cC^{[n]}\to|\beta|$ has fibers of the same dimension $n$.
\end{pro}
\begin{proof}
    We need to prove that for every curve $C\in|\beta|$, the Hilbert scheme $C^{[n]}$ has dimension $n$. This can be checked locally in the analytic topology (see, for example, \cite[\S 2.1, p. 288]{dCM00}). For every $Z\in C^{[n]}$, $Z$ is a zero-dimensional closed subscheme of length $n$ in $C$. Suppose that the underlying set of $Z$ consists of $m$ distinct points $x_1,\cdots,x_m\in C$ ($m\in\bZ_{>0}$). Choose for each $j\in\{1,\cdots,m\}$ an analytic open neighborhood $U_j$ of $x_j$ in $S$. We may assume these neighborhoods are disjoint and that they are isomorphic to analytic open subsets of $\bC^2$. Then $U_j':=U_j\cap C$ is an analytic open subset of a plane curve and 
    $$\left(\coprod_{j=1}^m U_j'\right)^{[n]}=\coprod_{n_1+\cdots+n_m=n}(U_1')^{[n_1]}\times\cdots\times (U_m')^{[n_m]}$$
    is an analytic open neighborhood of $Z$ in $C^{[n]}$. By \cite[Theorem 1.1]{Luan23}, $(U_j')^{[n_j]}$ has dimension $n_j$, hence the result follows.
\end{proof}

\section{Families of one-dimensional sheaves} 
\label{sec:pic}
Throughout this section, 
$\beta\cdot H$ and $\chi\in\bZ$ are coprime. 
We will construct families of sheaves in $M_{\beta,\chi}$ over curves and use these families as testing curves to study the Picard number of $M_{\beta,\chi}$. 
There are three types of families that we are going to construct.
\subsection{Moving support}\label{moving curve}
Suppose there is a pencil 
$$\bP^1\subset|\beta|$$ such that $\bP^1\cap U_0\not=\emptyset$, where $U_0\subset |\beta|$ is the locus of smooth curves. There is a multisection $R_0\hookrightarrow h_{\beta,\chi}^{-1}(\bP^1)$ of $h_1: h_{\beta,\chi}^{-1}(\bP^1)\rightarrow \bP^1$ of degree $d_0\in\bZ_{>0}$, i.e., the composition $h_1|_{R_0}:R_0\to\bP^1$ is a finite morphism of degree $d_0$. Let $R$ be the normalization of $R_0$ and let \[ \tau: R \rightarrow M_{\beta,\chi} \] be the composition of the normalization map and the closed embedding $R_0 \hookrightarrow M_{\beta,\chi}$.
Recall that $\lambda([\cO_x])=h_{\beta,\chi}^*\cO_{|\beta|}(1)$ by Corollary \ref{cor:lambda0}. Then we have
\begin{equation}
    \label{degmov}   \deg  (\tau^{*}\lambda([\cO_x]))=d_0>0.
\end{equation}

\subsection{Fixed smooth support}\label{smoothfiber}
We fix a smooth and connected curve $C_0\in\left|\beta\right|$. Let $\nu:C_0\times C_0\hookrightarrow C_0\times S$ be the inclusion map and let $y\in C_0$ be a point (also viewed as a point of $S$). We construct a family of sheaves by varying a point on $C_0$ while fixing $p_a(\beta)+\chi-2$ points as follows.
Let 
$$\cF_{C_0}:=\nu_*\cO_{C_0\times C_0}((p_a(\beta)+\chi-2)\,C_0\times\{y\}+\D_0),$$
where $\D_0$ is the diagonal in $C_0\times C_0$. Then under projection $C_0 \times S \rightarrow C_0$,  $\cF_{C_0}$ is a $C_0$-flat family of sheaves in $h_{\beta,\chi}^{-1}(C_0)$. This is a direct generalization of the construction in \cite[\S 2.2, p. 491]{CC15}.

To compute the Chern character $\mathrm{ch}(\cF_{C_0})$, we apply the Grothendieck--Riemann--Roch formula and obtain
\begin{equation*}
    \begin{aligned}
        &\mathrm{ch}(\cF_{C_0})\mathrm{td}(C_0\times S)\\
        =&\,\nu_*[\mathrm{ch}(\cO_{C_0\times C_0}((p_a(\beta)+\chi-2)\,C_0\times\{y\}+\D_0))\mathrm{td}(C_0\times C_0)]\\
        =&\, \nu_*[(1+(p_a(\beta)+\chi-2)\,C_0\times\{y\}+\D_0+(\chi-1)\{y\}\times\{y\})\mathrm{td}(C_0\times C_0)],
    \end{aligned}
\end{equation*}
where the second equality follows from $\D_0^2=(2-2p_a(\beta))\{y \}\times\{y \}$ (for simplicity of notation we do not distinguish different zero-cycles).
Therefore,
\begin{equation}
\label{grrc}
        \mathrm{ch}(\cF_{C_0})p_S^*\mathrm{td}(S)=C_0\times C_0+(\chi-1)C_0\times\{y\}+\nu_*(\D_0)+(\chi-p_a(\beta))\{y\}\times\{y\},
\end{equation}
where $p_S:C_0\times S\to S$ and $p_2:C_0\times C_0\to C_0$ are the projections to the second factors.

\subsection{Fixed reducible support}\label{singularfiber}
Suppose $D_0=C_1+C_2\in\left|\beta\right|$, where $C_1$ and $C_2$ are two smooth, connected curves on $S$. Assume 
that $C_1$ intersects $C_2$ transversally at $$n_0:=C_1\cdot C_2\geq |\chi|+1$$ points, say at $C_1\cap C_2=\{y_1,\cdots,y_{n_0}\}$. Then there is a natural exact sequence of sheaves on $D_0$
\begin{equation*}
    0 \lra \cO_{D_0} \lra \cO_{C_1} \bigoplus  \cO_{C_2} \lra \mathop{\bigoplus}_{j=1}^{n_0} \bC_{y_j} \lra 0.
\end{equation*}
Denote by $g_i$ ($i=1,2$) the genus of $C_i$.
Then by the adjunction formula,
$$p_a(\beta)=g_1+g_2+n_0-1.$$

Let $\cL_i$ be a line bundle on $C_i$ with $\chi(\cL_i)=\chi_i$, where
\begin{equation*}
\label{chiL12}
    \chi_1=\max\{\chi,0\}+1\quad\textnormal{and}\quad\chi_2=\chi+n_0+1-\chi_1.
\end{equation*}
We can glue $\cL_i$ to obtain a line bundle $\cL$ on the nodal curve $D_0$ with 
\begin{equation*}\label{chiL}
\chi(\cL)=\chi_1+\chi_2-n_0=\chi+1
\end{equation*} 
as follows. At $y_j$ ($j=1,\cdots,n_0$),
we fix the isomorphisms $\cL_1|_{y_j}\cong\cL_2|_{y_j}\cong \bC_{y_j}$. Let $\cL$ be the kernel of the following surjection induced by evaluation maps
$$\pi_{D_0}:\iota_{1*}\cL_1\bigoplus \iota_{2*}\cL_2 \lra \mathop{\bigoplus}_{j=1}^{n_0} \bC_{y_j},$$
where $\iota_i:C_i\hookrightarrow D_0$ is the inclusion map. By a local calculation, $\cL$ has constant rank $1$ at every point of the reduced curve $D_0$, hence $\cL$ is a line bundle on $D_0$ with $\cL|_{C_i}\cong \cL_i$ ({\it cf.} \cite[Proposition 3.2]{FB21}). In particular, $\iota_{0\ast} \cL$ is a one-dimensional sheaf on $S$, where $\iota_0:D_0 \hookrightarrow S$ is the closed embedding.

\begin{lem}
\label{lem:Lst}
    The one-dimensional sheaf $\iota_{0\ast} \cL$ is stable with respect to $H$.
\end{lem}
\begin{proof}
    Since $\cL|_{C_i}\cong\cL_i$ by the construction of $\cL$, we have short exact sequences 
    $$0\lra\cK_2\lra\cL\stackrel{r_1}\lra\iota_{1\ast}\cL_1\lra0\quad
   \text{and} \quad0\lra\cK_1\lra\cL\stackrel{r_2}\lra\iota_{2\ast} \cL_2\lra0.$$
    By the construction of $\cL$, the slopes of $\iota_{0\ast}\cK_i$ are
    \begin{equation*}
    \begin{aligned}
        \mu(\iota_{0\ast}\cK_2)&=\frac{\chi(\cL)-\chi(\cL_1)}{C_2\cdot H}=\frac{\min\{\chi,0\}}{C_2\cdot H},\\
        \mu(\iota_{0\ast}\cK_1)&=\frac{\chi(\cL)-\chi(\cL_2)}{C_1\cdot H}=\frac{\max\{\chi,0\}+1-n_0}{C_1\cdot H}.
    \end{aligned}
    \end{equation*}
    In particular, it follows from our assumption $n_0\geq|\chi|+1$ that
    \begin{equation}
    \label{eq:muki}
        \mu(\iota_{0\ast}\cK_i)<\frac{\chi}{\beta\cdot H}\quad(i=1,2).
    \end{equation}
    Let $\cF$ be a nonzero proper subsheaf of $\iota_{0\ast}\cL$. If $\Supp(\cF)=D_0$, then 
    $$\mu(\cF)=\frac{\chi(\cF)}{\beta\cdot H}<\frac{\chi(\cL)}{\beta\cdot H}=\mu(\iota_{0\ast}\cL).$$
    Now assume $\mathrm{Supp}(\cF)=C_1$ or $C_2$. Consider the restriction $r_i|_{\cF}:\cF\to\iota_{i\ast}\cL_i$, where $\cF$ is viewed as a subsheaf of $\cL$.
    If $r_1|_{\cF}$ is zero, then $\cF$ is a subsheaf of $\iota_{0\ast}\cK_2$. Since $\cK_2$ is a line bundle supported on $C_2$, by (\ref{eq:muki}) we have
    $$\mu(\cF)\leq\mu(\iota_{0\ast}\cK_2)<\mu(\iota_{0\ast}\cL)=\frac{\chi+1}{\beta\cdot H}.$$
    If $r_1|_{\cF}$ is nonzero, then $\mathrm{Supp}(\cF)=C_1$ and $r_2|_{\cF}=0$. Hence $\cF$ is a subsheaf of $\iota_{0\ast}\cK_1$ and again by (\ref{eq:muki}),
    $$\mu(\cF)\leq\mu(\iota_{0\ast}\cK_1)<\mu(\iota_{0\ast}\cL),$$
    which completes the proof.
\end{proof}

Next, we will construct a $C_1$-flat family of sheaves in $h_{\beta,\chi}^{-1}(D_0)$. The idea is to make modifications of $\cL$ constructed above along points on $C_1$. To do so, let $\cE_0:=p_S^*(\iota_{0*}\cL)$ be the pull back along the projection $p_S:C_1\times S\to S$. Then there is a short exact sequence of sheaves on $C_1\times S$
$$0\lra\cK\lra\cE_0\lra \nu_{1*}(\cE_0|_{C_1\times C_1})\lra0,$$
where $\nu_1:C_1\times C_1\to C_1\times S$ is the inclusion.
Note that $\cE_0|_{C_1\times C_1}\cong p_2^*\cL_1$
is a line bundle, where $p_2:C_1\times C_1\to C_1$ is the projection to the second factor.

Write 
$$\cG_0:=(\cE_0|_{C_1\times C_1})\otimes \cO_{C_1\times C_1}(-\D_{1})\subset\cE_0|_{C_1\times C_1},$$
where $\D_{1}$ is the diagonal in $C_1\times C_1$.
Let $\cF_{C_1}$ be the subsheaf of $\cE_0\oplus \nu_{1*}\cG_0$ defined by
$$\cF_{C_1}(V)=\{(s,t)\in\cE_0(V)\oplus \cG_0(\nu_1^{-1}(V)):s|_{\nu_1^{-1}(V)}=t\in\cE_0|_{C_1\times C_1}(\nu_1^{-1}(V))\}$$
for every open subset $V$ of $C_1\times S$.
Then
we have a commutative diagram of $C_1$-flat sheaves with exact rows
\begin{equation}
\label{diagram:ec1}
\begin{tikzcd}
{0}\arrow[r,""] & {\cK} \arrow[r,""] \arrow[d, "="]           & {\cF_{C_1}} \arrow[r,""] \arrow[d, hook, ""]  & {\nu_{1*}\cG_0}\arrow[d, hook, ""]\arrow[r,""] & {0}\\
{0}\arrow[r,""] &{\cK} \arrow[r, ""]  & {\cE_0} \arrow[r,""] & {\nu_{1*}(\cE_0|_{C_1\times C_1})} \arrow[r,""] & {0}.
\end{tikzcd}
\end{equation}

\begin{lem}
\label{lem:fc1}
    The sheaf $\cF_{C_1}$ is a $C_1$-family of sheaves in $h_{\beta,\chi}^{-1}(D_0)$.
\end{lem}
\begin{proof}
    For every point $x\in C_1$, restricting the diagram (\ref{diagram:ec1}) of $C_1$-flat sheaves to $\{x\}\times S\cong S$ yields
    \begin{equation*}
\begin{tikzcd}
{0}\arrow[r,""] & {\iota_{0\ast}\cK_2} \arrow[r,""] \arrow[d, "="]           & {\cF_{C_1}|_{\{x\}\times S}} \arrow[r,"r_1'"] \arrow[d, hook, ""]  & {\iota_{0\ast}\iota_{1\ast}\cL_1(-x)}\arrow[d, hook, ""]\arrow[r,""] & {0}\\
{0}\arrow[r,""] &{\iota_{0\ast}\cK_2} \arrow[r, ""]  & {\iota_{0\ast}\cL} \arrow[r,"r_1"] & {\iota_{0\ast}\iota_{1\ast}\cL_1} \arrow[r,""] & {0}.
\end{tikzcd}
\end{equation*}
As a subsheaf of $\iota_{0\ast}\cL$, the torsion sheaf $\cF_{C_1}|_{\{x\}\times S}$ is one-dimensional with determinant $\det(\cF_{C_1}|_{\{x\}\times S})=\det(\iota_{0\ast}\cL)\cong\cO_S(\beta)$ and Euler characteristic
$$\chi(\cF_{C_1}|_{\{x\}\times S})=\chi(\cK_2)+\chi(\cL_1(-x))=\chi(\cL)-1=\chi.$$
We claim that $\cF_{C_1}|_{\{x\}\times S}$ is stable. Let $\cF$ be a nonzero proper subsheaf of $\cF_{C_1}|_{\{x\}\times S}$. As in the proof of Lemma \ref{lem:Lst}, we may assume $\mathrm{Supp}(\cF)$ is $C_1$ or $C_2$.
Since $\cF$ is also a subsheaf of $\iota_{0\ast}\cL$, the proof of Lemma \ref{lem:Lst} implies that 
$$\mu(\cF)\leq \max\{\mu(\iota_{0\ast}\cK_1),\,\mu(\iota_{0\ast}\cK_2)\}<\frac{\chi}{\beta\cdot H}=\mu(\cF_{C_1}|_{\{x\}\times S})$$
by (\ref{eq:muki}). Hence $\cF_{C_1}|_{\{x\}\times S}$ is stable and the result follows.
\end{proof}
\begin{rem}
    Indeed, our construction of the family $\cF_{C_1}$ can  be viewed as a family of Hecke modifications of $\cL$ at points on $C_1$. 
\end{rem}
Using (\ref{diagram:ec1}) and the Grothendieck--Riemann--Roch formula, as in (\ref{grrc}) we have
\begin{equation}
\label{grrc1}
    \begin{aligned}
        &(\mathrm{ch}(\cF_{C_1})-\mathrm{ch}(\cE_0))p_S^*\mathrm{td}(S)\\
        =&\,(\mathrm{ch}(\nu_{1*}\cG_0)-\mathrm{ch}(\nu_{1*}(\cE_0|_{C_1\times C_1}))p_S^*\mathrm{td}(S)\\
        = &\,\nu_{1*}[p_2^*\mathrm{ch}(\cL_1)(\mathrm{ch}(\cO_{C_1\times C_1}(-\D_{1}))-1)p_2^*\mathrm{td}(C_1)]\\
        = &-\D_{1}+(1-g_1-\chi_1)\{y\}\times\{y\}.
    \end{aligned}
\end{equation}

\subsection{Intersection numbers}
\label{subsec:internum}
By the assertion 3 in \cite[Theorem 8.1.5]{HL10}, if $C$ is a curve and $\cF_C$ is a $C$-flat family of sheaves in $M_{\beta,\chi}$ inducing a morphism $\phi_{\cF_C}:C\to M_{\beta,\chi}$, then
$$\phi_{\cF_C}^*\lambda(\mathfrak{a})=\det(p_{C\,!}([\cF_C]\cdot p_S^*\mathfrak{a}))\quad (\mathfrak{a}\in K_{\beta,\chi}^0),$$
where $\lambda$ is as in (\ref{lambda}) and $p_C$ ({\it resp.} $p_S$) is the 
projection from $C\times S$ to $C$ ({\it resp.} $S$).
Thus the intersection number of $\lambda(\mathfrak{a})$ with $C$ is given by
\begin{equation}\label{deg}
\begin{aligned}
    &\deg(\phi_{\cF_C}^*\lambda(\mathfrak{a}))=\deg(\det(p_{C\,!}([\cF_C]\cdot p_S^*\mathfrak{a})))\\
    =&\deg( c_1((p_{C\,!}([\cF_C]\cdot p_S^*\mathfrak{a})))\\
    \stackrel{(\ast)}=&\deg (p_{C*}\{\mathrm{ch}(\cF_C)\cdot p_S^*\mathrm{ch}(\mathfrak{a})\cdot p_S^*\mathrm{td}(S)\}_3)\\
    =&\deg \left[ \{\mathrm{ch}(\cF_C)\cdot p_S^*\mathrm{td}(S)\}_1 \cdot \ch_2( p_S^*\mathfrak{a})+\right.\\
    & \left. \{\mathrm{ch}(\cF_C)\cdot p_S^*\mathrm{td}(S)\}_2\cdot \ch_1( p_S^*\mathfrak{a}) +r(\mathfrak{a})  \{\mathrm{ch}(\cF_C)\cdot p_S^*\mathrm{td}(S)\}_3\right] .
\end{aligned}
\end{equation}
where $\{\}_i$ means taking the degree $i$ part of a Chow class, $r(\mathfrak{a})$ is the rank of $\mathfrak{a}$ and $(\ast)$ follows from the Grothendieck--Riemann--Roch formula.

Now we use our construction of the testing curves and the calculations in \S \ref{moving curve}-\S \ref{singularfiber} to compute the intersection numbers of determinant line bundles with these testing curves. 
Let $x$ be a point of $S$ and let $L$ be a divisor on $S$. Then a direct calculation shows \[[\cO_x]\in K_{\beta,\chi}^0,\ \hbox{and}\  \theta_L:=[-(\beta\cdot L)\cO_S+\chi\cO_L]\in K_{\beta,\chi}^0.\]
Note that
\begin{equation}
\label{ch:kc}
    \mathrm{ch}([\cO_x])=[x] \quad\textnormal{and}\quad\mathrm{ch}(\theta_L)=-\beta\cdot L+\chi L-\frac{\chi}{2}L^2.
\end{equation}

\begin{enumerate}
    \item If $C$ is the curve $R$ with $\phi_{\cF_R}=\tau:R\to M_{\beta,\chi}$ as constructed in \S \ref{moving curve}, then  by (\ref{degmov}), 
$$\deg\phi_{\cF_R}^*\lambda([\cO_x])=d_0>0.$$
    \item If $C$ is the curve $C_0$ constructed in \S \ref{smoothfiber}, then 
by (\ref{grrc}), (\ref{deg}), and (\ref{ch:kc}),
\begin{equation*}
\begin{split}
&\deg\phi_{\cF_{C_0}}^*\lambda([\cO_x])=0,\\ 
&\deg\phi_{\cF_{C_0}}^*\lambda(\theta_L)=p_a(\beta)\beta\cdot L.
        \end{split}
    \end{equation*}

    \item Let $C$ be the curve $C_1$ constructed in \S \ref{singularfiber} (see Lemma \ref{lem:fc1}).
    As $\cE_0=p_S^*(\iota_{0*}\cL)$, for each $\mathfrak{a}\in K_{\beta,\chi}^0$, we have
    \begin{equation}
    \label{cE0}
        \{\mathrm{ch}(\cE_0)\cdot p_S^*\mathrm{ch}(\mathfrak{a})\cdot p_S^*\mathrm{td}(S)\}_3=\{p_{S}^*\mathrm{ch}(\iota_{0*}\cL)\cdot p_S^*\mathrm{ch}(\mathfrak{a})\cdot p_S^*\mathrm{td}(S)\}_3=0.
    \end{equation}
    Thus by (\ref{grrc1}), (\ref{deg}), (\ref{ch:kc}), and (\ref{cE0}),
    \begin{equation*}
        \begin{split}
           &\deg\phi_{\cF_{C_1}}^*\lambda([\cO_x])=0 , \\
           &\deg\phi_{\cF_{C_1}}^*\lambda(\theta_L)=-\chi C_1\cdot L+(g_1+\chi_1-1)\beta\cdot L. 
        \end{split}
    \end{equation*}
\end{enumerate}

\subsection{Lower bound of Picard numbers}

Recall that $[L_1],\cdots,[L_{\rho}]$ is a basis for the $\bQ$-vector space
$\mathrm{Num}_{\bQ}(S)=\Pic(S)_\bQ/\equiv_{\mathrm{num}}$, 
where $L_i\in\Pic(S)$ $(i=1,\cdots, \rho:=\rho(S))$ are smooth, connected effective divisors.
We use the following determinant line bundles
\begin{equation}
\label{def:lambda_i}
 \lambda_0:=\lambda([\cO_x]),\ \ \ \lambda_i:=\lambda(\theta_{L_i})
\end{equation}
to prove the lower bound, assuming $\beta$ satisfies $(\mathrm{P})$.

\begin{thm}\label{picardnum}
Suppose $\beta\cdot H$ and $\chi$ are coprime.
When $\beta$ satisfies $(\mathrm{P})$ as defined in Definition \ref{P}, we have the following relation of Picard numbers
    $$\rho(M_{\beta,\chi})\geq\rho(S)+1.$$
\end{thm}
 \begin{proof} 
We claim that $[\lambda_0],[\lambda_1],\cdots,[\lambda_\rho]$ which are defined in (\ref{def:lambda_i}) are linearly independent in $\mathrm{Num}_{\bQ}(M_{\beta,\chi})$, thus the result follows. 
We prove the claim by contradiction. If to the contrary that $[\lambda_0],[\lambda_1],\cdots,[\lambda_\rho]$ were linearly dependent, there would exist $a_0,a_1,\cdots,a_\rho\in\bQ$, not all zero, such that $$a_0\lambda_0+a_1\lambda_1+\cdots+ a_\rho\lambda_\rho\equiv_{\mathrm{num}}0.$$ 
Then for every curve $C$ and every morphism $\phi_{\cF_C}:C\to M_\beta$, 
\begin{equation} \label{singmatr}
    a_0\deg \phi_{\cF_C}^*\lambda_0+a_1\deg \phi_{\cF_C}^*\lambda_1+\cdots+a_\rho\deg \phi_{\cF_C}^*\lambda_\rho=0.
\end{equation}

Write $\beta\equiv_{\mathrm{num}}\sum_{i=1}^\rho d_iL_i$ with $d_i\in\bQ$. We may assume $d_1\not=0$.
Let $L_2',\cdots,L_r'$ be given as in Definition \ref{P} (3) and let $D_j=L_j+L_j'$ ($j=2,\cdots, \rho$) be the nodal curves that play the role of $D_0$ in \S \ref{singularfiber}.
By the calculations (1)-(3) at the end of \S \ref{subsec:internum}, we have the following table, where each entry is the intersection number (\ref{deg}) of a line bundle (a term in the first row) with a curve (a term in the first column), $n_j\in\bZ$ and the $*$'s are unimportant numbers. 

\begin{equation*}
\renewcommand\arraystretch{1.3}
    \begin{array}{c|cccc}
\cdot & \lambda_0  & \lambda_1&  \cdots &   \lambda_\rho \\ \hline 
R & d_0 &  * & \cdots  & *  \\
C_0 & 0& p_a(\beta)\beta\cdot L_1 & \cdots & p_a(\beta) \beta\cdot L_\rho \\
L_2 & 0&  -\chi L_2\cdot L_1+n_2\beta\cdot L_1  & \cdots  &  -\chi L_2\cdot L_\rho+n_2\beta\cdot L_\rho  \\
\vdots &   & &  &   \\
L_{\rho} & 0 &  -\chi L_\rho\cdot L_1+n_\rho\beta\cdot L_1  & \cdots  &  -\chi L_\rho\cdot L_\rho+n_\rho\beta\cdot L_\rho  \\
\end{array}
\end{equation*}

Thus the determinant of the intersection matrix can be calculated as

\[
\begin{aligned}
    &\det
\begin{pmatrix}
  d_0 & * & *   & *\\[5pt]
  0& p_a(\beta)\beta\cdot L_1  & \cdots & p_a(\beta) \beta\cdot L_\rho\\[5pt]
  0&  -\chi L_2\cdot L_1+n_2\beta\cdot L_1 & \cdots  &  -\chi L_2\cdot L_\rho+n_2\beta\cdot L_\rho\\[5pt]
  \vdots & \vdots & \vdots & \vdots\\[5pt]
  0 &  -\chi L_\rho\cdot L_1+n_\rho\beta\cdot L_1 & \cdots  &  -\chi L_\rho\cdot L_\rho+n_\rho\beta\cdot L_\rho
\end{pmatrix}\\[6pt]
=&(-\chi)^{\rho-1}d_0\cdot p_a(\beta)\det
\begin{pmatrix}
   \beta\cdot L_1 &  \cdots &  \beta\cdot L_\rho\\[5pt]
    L_2\cdot L_1 &\cdots  &  L_2\cdot L_\rho\\[5pt]
   \vdots & \vdots & \vdots \\[5pt]
   L_\rho\cdot L_1 & \cdots &  L_\rho\cdot L_\rho
\end{pmatrix}\\[6pt]
=&(-\chi)^{\rho-1}d_0d_1\cdot p_a(\beta)\det(L_i\cdot L_k)_{1\leq i,k\leq\rho},
\end{aligned}
\]
which is nonzero by the choice of $L_i$, contradicting (\ref{singmatr}).
 \end{proof}



\section{Wall-crossing approach to computing divisor class groups}
\label{sec:wc}

To avoid irreducibility issues as indicated in Section \ref{sec:other}, we denote by $M_{\beta,\chi}^+\subset M_{\beta,\chi}$ 
 the irreducible component in which a general point $F\in M_{\beta,\chi}^+$ is a stable one-dimensional sheaf with $\mathrm{Supp}(F)$ a smooth curve in $|\beta|$.
 
We propose to use the wall-crossing technique of $\delta$-stable pairs to study its Chow group $\CH^1(M_{\beta,\chi}^+)_\bQ$ with rational coefficients. This technique was used in \cite{CGKT20} to compute Poincaré polynomials of moduli spaces of one-dimensional sheaves on del Pezzo surfaces when $p_a(\beta)\leq2$. In some nice cases, we determine $\CH^1(M_{\beta,\chi})_\bQ$. 

Denote by $M_{\beta,\chi}^\circ$ the open subscheme parametrizing $F\in M_{\beta,\chi}$ with $\mathrm{Supp}(F)$ an integral curve. By Proposition \ref{prop:smsupp}, $M_{\beta,\chi}^\circ\subset M_{\beta,\chi}^+$.
 We summarize the main results of this section as follows. 
\begin{thm} \label{thm4}
Assume $\codim(M_{\beta,\chi}^+\setminus M_{\beta,\chi}^\circ,M_{\beta,\chi}^+)
\geq 2$ and $\chi>1$.
Let 
\begin{equation}
\label{eq:d}
    d=\frac{K_S\cdot\beta+\beta^2}{2}+\chi.
\end{equation}
Then we have
$$\CH^1(M_{\beta,\chi}^+)_{\bQ}\oplus\bQ=\CH^1(\cC^{[d]})_{\bQ}.$$
Furthermore, we have $h^0(\cO_S(\beta))\ge d+1$ and
$$\CH^1(M_{\beta,\chi}^+)_{\bQ}\oplus\bQ=\begin{cases}
  \CH^1(S^{[d]})_\bQ, & \text{($h^0(\cO_S(\beta))=d+1$),}  \\
  \CH^1(S^{[d]})_\bQ\oplus\bQ, & \text{($h^0(\cO_S(\beta))>d+1$),} \end{cases}$$
if either of the following holds:
\begin{enumerate}
     \item $\beta$ is $(d-1)$-very ample.

    \item $d=2$, and $\beta$ induces a finite, birational morphism $\phi:S\rightarrow |\beta|$ onto its image.
    \item $d=3$, and $\beta$ induces a finite, birational morphism $\phi:S\rightarrow |\beta|$ onto its image $\phi(S)$ that contains at most finitely many lines and lies in a quadric hypersurface in $|\beta|$. 
\end{enumerate}   
\end{thm}

The proof of Theorem \ref{thm4} will be provided at the end of \S \ref{sec:4.2.2}.

\subsection{Moduli of $\delta$-stable pairs}
First we recall the notion of pairs that originated from the work of Le Potier \cite{LPo93}. It is a special case of the $\alpha$-semistable coherent system in \cite{He98}.
\begin{defn}
Let $\delta\in\bQ_{>0}$. A pair $(F,s)$ is called 
$\delta$-semistable if 
\begin{enumerate}
    \item $F$ is a one-dimensional sheaf with a section $s\in H^0(F)\setminus\{0\}$.
    \item For every nonzero proper subsheaf $G$ of $F$, 
    \begin{equation}
    \label{deltastab}
        \frac{\chi(G)+\epsilon(s,G)\delta}{c_1(G)\cdot H} \leq \frac{\chi(F)+\delta}{c_1(F)\cdot H},
    \end{equation}
    where 
    $$\epsilon(s,G)=\begin{cases}
    1,&\text{if $s:\cO_S\to F$ factors through $G$;}\\
    0,&\text{otherwise.}
    \end{cases}$$
\end{enumerate}
The above pair $(F,s)$ is called $\delta$-stable if it is $\delta$-semistable and the inequality (\ref{deltastab}) is strict. A $\delta$-semistable pair is called $\delta$-polystable if it is the direct sum of $\delta$-stable pairs.
\end{defn}    
Denote by $P_{\beta,\chi}(\delta)$ the coarse moduli space parametrizing $\delta$-polystable pairs $(F,s)$ with $\det F\cong\cO_S(\beta)$ and $\chi(F)=\chi$. As $\delta\in \bQ_{>0}$ varies, a wall-crossing phenomenon appears: there are finitely many rational numbers $0<w_1<\cdots < w_\ell<\infty$ such that for any $i\in\{0,1,\cdots,\ell\}$, and any $\delta,\delta'\in(w_i,w_{i+1})\cap\bQ$ ($w_0:=0,\,w_{\ell+1}:=\infty$), the moduli space $P_{\beta,\chi}(\delta)$ is isomorphic to $P_{\beta,\chi}(\delta')$. We denote \[ P_{\beta,\chi}(0^+):=P_{\beta,\chi}(\delta_0),\quad P_{\beta,\chi}(\infty):=P_{\beta,\chi}(\delta_\infty), \]
for $\delta_0\in(0,w_1)\cap\bQ$ and $\delta_\infty\in(w_\ell,\infty)\cap\bQ$.

When $\delta$ is sufficiently large, the condition of a pair $(F,s)$  being $\delta$-semistable is equivalent to that the cokernel of $s:\cO_S\to F$ is $0$-dimensional. In particular, $P_{\beta,\chi}(\infty)$ is the moduli space of Pandharipande--Thomas stable pairs \cite{PT09}. Then $P_{\beta,\chi}(\infty)$ is isomorphic to the relative Hilbert scheme $\cC^{[d]}$ of $d$ points (see \S\ref{subsec:relHil}) by \cite[Proposition B.8]{PT10}. It has the following property.

\begin{pro}
\label{prop:relhilirr}
If $d\geq0$, then the scheme $P_{\beta,\chi}(\infty)$ is integral of dimension
   $$h^0(\cO_S(\beta))-1+d,$$
with at worst locally complete intersection singularities.
\end{pro}
\begin{proof}
It suffices to prove the statement for $\cC^{[d]}$. 
Denote by $\cC_{U_0}^{[d]}$ the relative Hilbert scheme of points over the open subset $U_0\subset|\beta|$ parametrizing smooth curves. It is irreducible since fibers of $\pi^{[d]}$ over $U_0$ are irreducible of dimension $d$. It is smooth by Proposition \ref{prop:smsupp} and \cite[Proposition 14]{She12}, therefore contained in a single irreducible component $M$ of dimension $h^0(\cO_S(\beta))-1+d$. 
The proof of irreducibility of $P_{\beta,\chi}(\infty)$ is similar to that of \cite[Proposition 2.3]{PSSZ24}. If there is another irreducible component $M'$ other than $M$, then by Proposition \ref{eqdim:cC},
\[\dim M' <\dim M= h^0(\cO_S(\beta))-1+d. \]
On the other hand, the universal curve $\cC \subset S \times |\beta|$ induces a natural closed embedding $\cC^{[d]} \subset S^{[d]}\times |\beta|$. It follows from the argument of \cite[Section A.2]{KT14} that $\cC^{[d]}$ is the zero locus of a section of a rank $d$ vector bundle on $S^{[d]}\times |\beta|$, which implies 
\[ \dim M' \ge \dim(S^{[d]}\times|\beta|)-d=  h^0(\cO_S(\beta))-1+d, \]
contradiction! Hence $P_{\beta,\chi}(\infty)=M$ is irreducible.  And it is a local complete intersection in $S^{[d]}\times |\beta|$ since the expected dimension $d+ h^0(\cO_S(\beta))-1$ is achieved. In particular, $P_{\beta,\chi}(\infty)$ is generically smooth and Cohen--Macaulay, thus it is reduced.
\end{proof}


Denote by $P_{\beta,\chi}(\delta)^\circ$ the open subscheme parametrizing $(F,s)\in P_{\beta,\chi}(\delta)$ with $\mathrm{Supp}(F)$ an integral curve. We make a simple but useful observation about the wall-crossing. 

\begin{lem}
\label{lem:int}
    If $(F,s) \in P_{\beta, \chi}(c)^\circ$ for some $c\in\bQ_{>0}$, then $(F,s)$ is $\delta$-stable for any $\delta\in\bQ_{>0}$. In particular, $P_{\beta,\chi}(0^+)^\circ$ is isomorphic to $P_{\beta,\chi}(\infty)^\circ$.
\end{lem}
\begin{proof}
If  $(F,s)$ is not $\delta$-stable for some $\delta\in\bQ_{>0}$, then there is a proper destabilizing subsheaf $G\subset F$, i.e., 
\begin{equation*}
    \frac{\chi(G)+\epsilon(s,G)\delta}{c_1(G)\cdot H}\geq \frac{\chi(F)+\delta}{c_1(F)\cdot H}=\frac{\chi+\delta}{\beta\cdot H}.
\end{equation*}
By the purity of $F$ and the integrality of $\Supp(F)$, we have 
$$[\Supp(G)]=c_1(G)=c_1(F)=[\Supp(F)],$$ 
which implies $\chi(G)\geq\chi(F)$. But this leads to a contradiction since $F/G$ is a nonzero sheaf of dimension zero. 
\end{proof}

\subsection{Dimension estimates for bad loci}

\subsubsection{Estimates for large $\delta$}
Note that under the identification $P_{\beta,\chi}(\infty) \cong \cC^{[d]}$, we have
$$P_{\beta,\chi}(\infty)^\circ =(\pi^{[d]})^{-1}(U),$$ 
where $\pi^{[d]}:\cC^{[d]} \rightarrow |\beta|$ is as in \S \ref{subsec:relHil}. Since the morphism $\pi^{[d]}$ has fibers of the same dimension $d$ by Proposition \ref{eqdim:cC}, we have
\begin{equation}
\label{eq:codimp}
    \codim (P_{\beta,\chi}(\infty)\setminus P_{\beta,\chi}(\infty)^\circ,P_{\beta,\chi}(\infty))=\codim (|\beta|\setminus U,|\beta|).
\end{equation}  

Let $N:=h^0(\cO_S(\beta))-1$. Denote by $\pi_1: S^{[d]} \times S \rightarrow S^{[d]} $ and $\pi_2: S^{[d]} \times S \rightarrow  S$ the two projections. Then there is a natural morphism of vector bundles on $S^{[d]}$
\begin{equation*}
    \varphi:\ H^0(\cO_S(\beta))\otimes \cO_{S^{[d]}} \lra \cE:=\pi_{1\ast }(\pi^\ast_2 \cO_S(\beta) \otimes \cO_{\cZ}),
\end{equation*}
where $\cZ \subset S^{[d]} \times S$ is the universal closed subscheme.  Then the $k$-th degeneracy locus \[W^d_{\le k}(\beta):=\{Z\in S^{[d]}\ |\ \rank(\varphi|_Z) \le k\  \} \subset S^{[d]}\] is a closed subscheme with natural inclusions
\begin{equation*}
       \cdots \subset W^d_{\le k}(\beta)\subset W^d_{\le k+1}(\beta) \subset \cdots \subset W^d_{\le d-1}(\beta) \subset  W^d_{\le d}(\beta)=S^{[d]}.
\end{equation*}
By \cite[Theorem 14.4 (b)]{Ful98}, $W^{d}_{\le k}(\beta)$ has expected codimension $(N+1-k)(d-k)$. Write $W^d_k(\beta):=W^d_{\le k}(\beta)\setminus W^d_{\le k-1}(\beta)$. Then $ W^d_d(\beta)\subset S^{[d]}$ is an open subset. 


The fiber $\varphi|_Z$ of the bundle morphism $\varphi$ at a point $Z\in S^{[d]}$ is just the evaluation map 
\begin{equation}\label{evm}
  r_Z:\bC^{N+1}\cong H^0(\cO_S(\beta))\longrightarrow \bC^{d}\cong H^0(\cO_S(\beta)|_Z)
\end{equation}
induced by the short exact sequence for the closed subscheme $Z$ in $S$
\begin{equation}\label{ses}
    0 \lra \cI_Z\otimes \cO_S(\beta) \lra \cO_S(\beta) \lra \cO_S(\beta)|_Z \lra 0,
\end{equation}
where $\cI_Z$ is the ideal sheaf of $Z\subset S$.
Note that the morphism $\sigma_{d-1}:\cC^{[d]} \rightarrow S^{[d]}$ defined in (\ref{eq:defsigma}) has a natural $\bP^{N-k}$-bundle structure over $W^d_k$ (see the proof of Proposition \ref{kveryamp}). Our observation is that if 
\begin{equation}
\label{codimcond}
    \codim (W^d_{\le d-n}(\beta), S^{[d]} )\ge n+2\textnormal{ for every } n\in\bZ_{>0},
\end{equation}
then  $$\codim (\sigma_{d-1}^{-1}(W^d_{\le d-1}(\beta)),\cC^{[d]})\ge 2$$ and thus 
\begin{equation}
\label{eq:ch1cd}
\begin{aligned}
    \CH^1(\cC^{[d]})_\bQ&=\CH^1(\cC^{[d]}\setminus\sigma_{d-1}^{-1}(W^d_{\le d-1}(\beta)))_\bQ\\
    &=\begin{cases}
  \CH^1(W^d_{d}(\beta))_\bQ & \text{($N=d$)}  \\
  \CH^1(W^d_{d}(\beta))_\bQ\oplus\bQ & \text{($N>d$)} \end{cases}\\
    &=\begin{cases}
  \CH^1(S^{[d]})_\bQ & \text{($N=d$)}  \\
  \CH^1(S^{[d]})_\bQ\oplus\bQ & \text{($N>d$).} \end{cases}
\end{aligned}
\end{equation}
where the second identity follows from \cite[Theorem 3.3 (b)]{Ful98}.

Now we study the geometry of $W^d_{\le d-n}(\beta)$
from the perspective of the morphism \[\phi:S\to|\beta|=\bP(H^0(\cO_S(\beta))^\vee)\] induced by $\beta$. Assume $Z\in W^d_{\le d-1}(\beta)$ is the union of $d$ distinct points $x_1,\cdots,x_d$ on $S$. Then the evaluation map (\ref{evm}) has rank $\le d-1$ if and only if the $d$ points $\phi(x_1),\cdots,\phi(x_d)$ lie in a $(d-2)$-dimensional linear subspace of $|\beta|$.

For a partition $\mu=(n_1,\cdots,n_m)$ of $d$ (i.e., $m,\,n_i\in\bZ_{>0}$, $n_1\geq\cdots\geq n_m$ such that $\sum_{i=1}^mn_i=d$), denote by $S_\mu^{[d]}\subset S^{[d]}$ the subscheme that parametrizes all the closed subschemes corresponding to $0$-cycles of the form $\sum_{i=1}^mn_i\cdot z_{i}$, where $z_{i}$ are distinct points on $S$. 

\begin{pro}[length $2$]
\label{prop:lth2}
    Suppose the morphism $\phi:S\to|\beta|$ induced by $\beta$ is birational and finite onto its image. Then (\ref{codimcond}) holds for $d=2$.
\end{pro}
\begin{proof}
Since $\beta$ is base-point-free, $W^d_{\le d-2}(\beta)=W^2_{\le 0}(\beta)$ is empty by the description of (\ref{evm}). It remains to show $\dim W^{2}_{\le 1}(\beta)\le 1$.
Note that $S^{[2]}=S^{[2]}_{(1,1)} \sqcup S^{[2]}_{(2)}$. 
A union of two distinct points $x_1,\,x_2$ on $S$ belongs to $W^2_{\leq 1}(\beta)$ if and only if $\phi(x_1)=\phi(x_2)$. By our assumption on $\phi$, $$\dim (W^{2}_{\le 1}(\beta) \cap S^{[2]}_{(1,1)})\leq 1.$$

To estimate $\dim(W^{2}_{\le 1}(\beta) \cap S^{[2]}_{(2)})$, we take a closer look at the evaluation map (\ref{evm}) when $Z\in S^{[2]}_{(2)}$. Note that to give $Z\in S^{[2]}_{(2)}$ is equivalent to giving a point $x\in S$ (the underlying set of $Z$) and a nonzero tangent vector $v\in T_{x}S:=(m_x/m_x^2)^\vee$ up to a nonzero scalar (the length $2$ scheme structure) ({\it cf}. \cite[Ch. \uppercase\expandafter{\romannumeral2}, Exercise 2.8, p. 80]{Hartshorne}). We choose a basis $s_0,\cdots,s_N$ for $H^0(\cO_S(\beta))$ such that $s_0(x)\not=0$ and $s_1(x)=\cdots=s_N(x)=0$. Then the rank of (\ref{evm}) equals the rank of
$$\begin{pmatrix}
  1  & 0 & \cdots & 0 \\
  0 & v([s_1]) & \cdots & v([s_N])
\end{pmatrix},$$
where $[s_i]$ is viewed as an element in $m_x/m_x^2$ via a trivialization of $\cO_S(\beta)$ near $x$. Hence $Z=(x,v)\in W^{2}_{\leq 1}(\beta)$ if and only if $v$ is in the kernel of the tangent map at $x$ $$T_x\phi:T_xS\to T_{\phi(x)}|\beta|.$$
Since $\phi$ is finite and birational, we have
$$\dim\{z\in S \ |\ \rk T_z\phi\leq j\}\leq j\quad(j=0,1),$$
thus $\dim W^2_{\le 1}(\beta) \cap S^{[2]}_{(2)}\leq 1$.
\end{proof}
\begin{rem}
   The above assumption on $\beta$ to ensure $\codim(W^2_{\le 1}(\beta), S^{[2]})\ge 3$ is almost sharp. For example, if we only require that $\phi$ is birational onto its image, then $\phi$ may contract a curve $C$ to a point. In this case, the subscheme of $C^{[2]}$ that parametrizes two distinct points on $C$ is contained in $W^2_{\le 1}(\beta)$ and thus $codim (W^2_{\le 1}(\beta),S^{[2]}) \le 2$. 
\end{rem}
\begin{pro}[length $3$]
\label{prop:lth3}
   Assume $\phi:S\to|\beta|$ induced by $\beta$ is birational and finite onto its image. If the image $S':=\phi(S)$ contains at most finitely many lines in $|\beta|$ and lies in a quadric hypersurface, then (\ref{codimcond}) holds for $d=3$.
\end{pro}
\begin{proof}
We begin by showing $\dim W^3_{\le1}(\beta)\leq 2$. Consider the incidence variety $S^{[2,3]}$ of all pairs $(\xi,\xi')\in S^{[2]}\times S^{[3]}$ satisfying $\xi\subset \xi'$. Let $\rho_2:S^{[2,3]}\to S^{[2]}$ and $\rho_3:S^{[2,3]}\to S^{[3]}$ be the natural projections. Then $\rho_2\rho_3^{-1} W_{\le1}^3(\beta)\subset W^2_{\le1}(\beta)$. Note that $\rho_2^{-1}(\xi)$ is the blow-up of $S$ along $\xi$, so $\dim(\rho_2^{-1}(\xi)\cap\rho_3^{-1}W_{\le1}^3(\beta))\le1$ for a point $\xi\in W^2_{\le1}(\beta)$ since by adding a general point of $S$ to $\xi\in W^2_{\le1}(\beta)$ we get a point in $W^3_{2}(\beta)$. Therefore, by Proposition \ref{prop:lth2} we have
$$\dim W^3_{\le1}(\beta)-1\le\dim\rho_2\rho_3^{-1} W_{\le1}^3(\beta)\le\dim W^2_{\le1}(\beta)\le1.$$

To study $W^3_{\leq2}(\beta)$, we deal with the first two strata of $S^{[3]}=S^{[3]}_{(1,1,1)}\sqcup S^{[3]}_{(2,1)}\sqcup S^{[3]}_{(3)}$ using the following commutative diagram
\begin{equation*}
    \begin{tikzcd}
        \Gamma \arrow[d, "q_1"]  \arrow[dr, "q_2"] & \\
         (S^2\setminus\D_S)\times S \arrow[r, dashed, "\sigma"] &  S^{[3]},
    \end{tikzcd}
\end{equation*}
where $\D_S\subset S^2$ is the diagonal, $\sigma$ is induced by the quotient map from $S^3$ to the symmetric product $S^{(3)}$, and $q_1$ is the blow-up along the big diagonal $$\{(x,y,z)\in S^3 \ | \ x\not=y,\text{ either }z=x\text{ or }z=y \}\subset(S^2\setminus\D_S)\times S.$$ The image of $q_2$ is $S^{[3]}_{(1,1,1)}\cup S^{[3]}_{(2,1)}$. Recall that for $Z=x_1+x_2+x_3\in S^{[3]}_{(1,1,1)}$, $Z\in W^3_{\leq2}(\beta)$ if and only if $\phi(x_1),\,\phi(x_2),\,\phi(x_3)$ lie in the same line. By the proof of Proposition \ref{prop:lth2}, for $Z'=(x,v)\cup \{y\}\in S^{[3]}_{(2,1)}$, where $x,\,y$ are distinct points on $S$ and $v\in T_{x}S\setminus\{0\}$, if we choose $s_0,\cdots,s_N$ as in Proposition \ref{prop:lth2}, then the rank of $\varphi|_{Z'}$ equals the rank of
$$\begin{pmatrix}
  1  & 0 & \cdots & 0 \\
  0 & v([s_1]) & \cdots & v([s_N])\\
  s_0(y)  & s_1(y) & \cdots & s_N(y)
\end{pmatrix}.$$
In particular, if $Z'\in W^3_{\leq 2}(\beta)$ and $v\notin\ker T_x\phi$, then either $\phi(x)=\phi(y)$, or the line joining $\phi(x)\not=\phi(y)$ lies in the projective tangent space of $S'$ at $\phi(x)$.

We next prove that $W:= p_{12}(q_1q_2^{-1}W^3_{\leq 2}(\beta))$ is a proper closed subset in $S^2\setminus\D_S$, where $p_{12}:(S^2\setminus\D_S)\times S\to S^2\setminus\D_S$ is the projection. Since $p_{12}$ and $q_1$ are proper, $W$ is closed in $S^2\setminus\D_S$. By the assumption on $S'$, there exist two distinct points $x_1,\,x_2\in S$ in the locus where $\phi$ is an isomorphism onto its image, such that the line $\ell_{12}$ joining $\phi(x_1)\not=\phi(x_2)$ is not contained in $S'$. Since $S'$ lies in a quadric hypersurface, $\ell_{12}$ meets $S'$ transversally and therefore $(x_1,x_2)\in(S^2\setminus\D_S)\setminus W$. Thus $\dim W\leq 3$. Note that the fiber dimension of $p_{12}\circ q_1|_{q_2^{-1}W_{\leq 2}(\beta)}$ over $(y_1,y_2)\in S^2\setminus\D_S$ is $2$ if and only if $\phi(y_1)=\phi(y_2)$. The fiber dimension is $1$ if and only if either
\begin{enumerate}
    \item the line joining $\phi(y_1)\not=\phi(y_2)$ is contained in $S'$, or
    \item the line joining $\phi(y_1)\not=\phi(y_2)$ is not contained in $S'$, but $T_{y_1}\phi=0$ or $T_{y_2}\phi=0$.
\end{enumerate}
Combining with Proposition \ref{prop:lth2} and our assumption on $S'$, we conclude that $$\dim W^3_{\leq 2}(\beta)\cap(S^{[3]}_{(1,1,1)}\cup S^{[3]}_{(2,1)})\leq 3.$$

The third stratum $S^{[3]}_{(3)}$ is an irreducible closed subscheme of $S^{[3]}$ of dimension $4$. To show $\dim W^3_{\leq2}(\beta)\cap S^{[3]}_{(3)}\leq 3$, it suffices to show $S^{[3]}_{(3)}\setminus W^3_{\leq 2}(\beta)\not=\emptyset$. Let $x_0\in S$ be a point such that $T_{x_0}\phi$ is injective. Then $\cI_{0}^2$, where $\cI_{0}$ is the ideal sheaf of $x_0$, defines a length $3$ closed subscheme $Z_0\in S^{[3]}_{(3)}$ with an exact sequence
$$0\lra \cI_{0}/\cI^2_{0}\lra\cO_{Z_0}\lra\cO_{x_0}\lra0.$$
If we identify $H^0(\cO_S(\beta)|_{Z_0})$ with $H^0(\cO_S(\beta)|_{x_0})\oplus H^0(\cO_S(\beta)\otimes\cI_{0}/\cI^2_{0})$ via the above exact sequence, then the injectivity of $T_{x_0}\phi$ implies $Z_0\not\in W^3_{\leq2}(\beta)$, so we are done.
\end{proof}

\begin{rem}
Take $\beta=K_S+4\beta'+\gamma$ for any very ample divisor $\beta'$ and nef divisor $\gamma$. By \cite[Theorem 1.8.60]{Lazars04}, $\beta$ satisfies Property $(N_1)$, which implies that defining equations of $\phi(S)$ are all quadratics.
\end{rem}

\subsubsection{Controlling the wall-crossing for small $\delta$ }\label{sec:4.2.2}
Now assume $\chi>0$. In what follows, we use the duality between $M_{\beta,\chi}$ and $M_{\beta,-\chi}$ to estimate the dimensions. By $\pm\chi$, we mean ``$\chi$ ({\it resp.} $-\chi$)''.
Suppose $\beta\cdot H$ and $\chi$ are coprime, so that there is a universal sheaf $\cE_{\pm}$ on $M_{\beta,\pm\chi}\times S$. 

The natural morphism defined by forgetting the section in each pair
\[f_{\pm\chi}:\  P_{\beta,\pm \chi}(0^+)^\circ
\  \rightarrow M_{\beta,\pm \chi}^\circ 
\]
is a $\bP^k$-bundle over the stratum
\begin{equation*}
    Z_{\pm \chi,k}:=\{ F\in M_{\beta,\pm \chi}^\circ
    \ |\ h^0(F)=k+1 \}
\end{equation*}
({\it cf}. \cite[Proposition 3.17]{CGKT20}). Viewing $P_{\beta,\pm \chi}(0^+)^\circ$ and $M_{\beta,\pm \chi}^\circ$ as the relative Hilbert scheme and the relative compactified Jacobian, respectively, one sees that $f_{\pm\chi}$ is in fact the Abel--Jacobi map. By the upper semicontinuity \cite[Ch. \uppercase\expandafter{\romannumeral3}, Theorem 12.8]{Hartshorne}, the subset 
\begin{equation*}
    \{\ [F]\in M_{\beta,\pm \chi}^\circ
    \ |\ h^0(\cE_{\pm}|_{\{[F]\}\times S})=h^0(F) \ge m\ \}
\end{equation*}
is closed for all $m\in \bZ_{\ge 0}$, thus $Z_{\pm\chi,k}$ is locally closed
and $Z_{\chi,k_{min}}$ is an open subset of $M_{\beta,\chi}^\circ$, where $$k_{min}:=\min\{ h^0(F)-1\ |\ F\in M_{\beta, \chi}^\circ
\}.$$ 
\begin{pro}
\label{prop:dim}
  We have  $k_{min}=\chi-1$ and 
  \begin{equation*}
      \codim (M_{\beta,\chi}^\circ\setminus Z_{\chi,\chi-1},M_{\beta,\chi}^\circ) \ge \chi+1.
  \end{equation*}
  Moreover, $M_{\beta,\chi}^\circ\setminus Z_{\chi,\chi-1}=\emptyset$ if $\chi> p_a(\beta)-1=\beta(\beta+K_S)/2$. If $\chi\leq p_a(\beta)-1$, then
  \begin{equation}
  \label{eq:codim2}
      \codim(f_\chi^{-1}(M_{\beta,\chi}^\circ\setminus Z_{\chi,\chi-1}),P_{\beta,\chi}(0^+)^\circ)= \chi.
  \end{equation}

\end{pro}
\begin{proof}
   Let $\iota:S\hookrightarrow \bP^n$ be the closed embedding induced by $mH$ for a large $m$. Then $M_{\beta,\chi}$ can be viewed as a subscheme of the moduli space of one-dimensional sheaves on $\bP^n$. By \cite[Theorem 13]{Maican10}, there is an isomorphism 
   \begin{equation}\label{dualiso}
      (-)^D: M_{\beta,\chi} \rightarrow M_{\beta,-\chi}
   \end{equation}
   by sending $F$ to $F^D:=\cext^1(F,\omega_S)$, which is compatible with the Hilbert--Chow morphisms. Note that $\cext^{n-1}(\iota_*F,\omega_{\bP^n})\cong\iota_*\cext^1(F,\omega_S)$ by Grothendieck--Verdier duality, so our notation $F^D$ is consistent with that in \cite{Maican10}. Moreover, by 
   \cite[Corollary 6]{Maican10},
   \begin{equation*}
       h^1(F^D)=h^0(F)= h^0(F^D)+\chi.
   \end{equation*}
   Thus the isomorphism (\ref{dualiso}) induces isomorphisms of the strata 
 $Z_{\chi,k} \cong Z_{-\chi,k-\chi}$
 and
 \begin{equation}
    \label{eq:codimz}
        \im(f_{-\chi})=\{E\in M_{\beta,-\chi}^\circ \ | \ h^0(E)\geq1\}\cong\{F\in M_{\beta,\chi}^\circ\ | \ h^0(F)\geq \chi+1\}.
    \end{equation}
If $\chi\geq p_a(\beta)$, then $P_{\beta,-\chi}(0^+)^\circ=P_{\beta,-\chi}(\infty)^\circ$ is empty ({\it cf}. \cite[(2.5)]{PT10}), which implies $\im(f_{-\chi})=\emptyset$ and $M_{\beta,\chi}^\circ=Z_{\chi,\chi-1}$. 

Now assume $\chi\leq p_a(\beta)-1$.
  By Propositions \ref{prop:smsupp}, \ref{prop:relhilirr}, and Lemma \ref{lem:int}, $$\dim P_{\beta,-\chi}(0^+)^\circ=\dim P_{\beta,-\chi}(\infty)^\circ=\dim|\beta|+(p_a(\beta)-1-\chi)=\dim M_{\beta,-\chi}^\circ-1-\chi.$$ Hence $f_{-\chi}$ is not surjective and the maximal open strata in $M_{\beta,-\chi}$ is $Z_{-\chi,-1}$, which implies $k_{\mathrm{min}}=\chi-1$.
    Combined with (\ref{eq:codimz}), it follows that $$\codim (M_{\beta,\chi}^\circ\setminus Z_{\chi,\chi-1},M_{\beta,\chi}^\circ)=\codim(\im(f_{-\chi}),M_{\beta,-\chi}^\circ)\geq \chi+1.$$
    Since $P_{\beta,\chi}(0^+)^\circ$ is irreducible, so is $\im (f_{-\chi})$. For every stratum $Z_{\chi,k}$ with $k\geq\chi$, 
    \begin{equation}
    \label{eq:zkn}
        \begin{aligned}
            \dim f_\chi^{-1}(Z_{\chi,k})&=\dim Z_{\chi,k}+k=(\dim Z_{-\chi,k-\chi}+k-\chi)+\chi\\
            &=\dim f_{-\chi}^{-1}(Z_{-\chi,k-\chi})+\chi\\
            &\leq\dim P_{\beta,-\chi}(0^+)^\circ+\chi=\dim P_{\beta,\chi}(0^+)^\circ-\chi,
        \end{aligned}
    \end{equation}
    where the last equality follows by Proposition \ref{prop:relhilirr} and Lemma \ref{lem:int}.
    If we take $\ell_0:=\min\{\ell\in\bZ_{\geq0}\ | \ Z_{-\chi,\ell}\not=\emptyset\}$, then the inequality in (\ref{eq:zkn}) is actually an equality for $k=\ell_0+\chi$ since $Z_{-\chi,\ell_0}$ is a dense open subset of $\im(f_{-\chi})$, which proves (\ref{eq:codim2}).
\end{proof}

Putting these all together, we obtain the following.
\begin{thm}
\label{thm:main41}
Assume $\codim(M_{\beta,\chi}^\circ,M_{\beta,\chi}^+)
\geq 2$ and $\chi>1$. Let $d$ be given as in (\ref{eq:d}). Then
$$\CH^1(M_{\beta,\chi}^+)_{\bQ}\oplus\bQ=\CH^1(\cC^{[d]})_{\bQ}.$$
If furthermore $\beta$ is ($d-1$)-very ample, then \footnote{If $h^{\geq1}(\cO_S(\beta))=0$, then $h^{0}(\cO_S(\beta))>d+1\iff\chi(\cO_S)>K_S\cdot\beta+\chi+1$ by Riemann--Roch.}
$$\CH^1(M_{\beta,\chi}^+)_{\bQ}\oplus\bQ=\begin{cases}
  \CH^1(S^{[d]})_\bQ, & \text{if $h^0(\cO_S(\beta))=d+1$,}  \\
  \CH^1(S^{[d]})_\bQ\oplus\bQ, & \text{if $h^0(\cO_S(\beta))>d+1$.} 
\end{cases}$$

\end{thm}
\begin{proof}
By Proposition \ref{prop:dim} and our assumption $\codim(M_{\beta,\chi}^\circ,M_{\beta,\chi}^+)
\geq 2$, we have 
\begin{equation}
\label{eq:pd1}
\CH^1(M_{\beta,\chi}^+)_{\bQ}=\CH^1(M_{\beta,\chi}^\circ)_{\bQ}=\CH^1(Z_{\chi,\chi-1})_\bQ.
\end{equation}
Let $p_Z: Z_{\chi,\chi-1} 
\times S \rightarrow Z_{\chi,\chi-1}
$ be the projection. By the cohomology and base change theorem \cite[Ch. III, Theorem 12.11]{Hartshorne}, the direct image $\cH:=p_{Z*}(\cE|_{Z_{\chi,\chi-1}\times S})$ of the universal sheaf $\cE$ on $M_{\beta,\chi}\times S$ is a vector bundle on $Z_{\chi,\chi-1}$ of rank $\chi$. 
Then by the modular interpretation, $f_{\chi}^{-1}(Z_{\chi,\chi-1})\cong\bP(\cH^\vee)$, and thus
\begin{equation}
\label{eq:pd2}
\begin{aligned}
    \CH^1(Z_{\chi,\chi-1})_\bQ\oplus\bQ&= \CH^1(f_{\chi}^{-1}(Z_{\chi,\chi-1}))_\bQ&(\text{by \cite[Theorem 3.3 (b)]{Ful98}})\\
    &=\CH^1(P_{\beta,\chi}(0^+)^\circ)_\bQ&(\text{by Proposition \ref{prop:dim}})\\
    &=\CH^1(P_{\beta,\chi}(\infty)^\circ)_\bQ & (\text{by Lemma \ref{lem:int}}).
\end{aligned}
\end{equation}
The upper semicontinuity of $h_{\beta,\chi}$ implies
$$\codim(|\beta|\setminus U,|\beta|)\geq\codim(M_{\beta,\chi}^\circ,M_{\beta,\chi}^+)\geq 2,$$
and therefore, by (\ref{eq:codimp}) we have
\begin{equation}
\label{eq:pd3}
    \CH^1(P_{\beta,\chi}(\infty)^\circ)_\bQ=\CH^1(P_{\beta,\chi}(\infty))_\bQ=\CH^1(\cC^{[d]})_\bQ.
\end{equation}
The first statement follows from (\ref{eq:pd1}), (\ref{eq:pd2}), and (\ref{eq:pd3}). 

If $\beta$ is $(d-1)$-very ample, $\cC^{[d]}$ is a projective bundle over $S^{[d]}$ with fibers of dimension $h^0(\cO_S(\beta))-d-1\geq0$ by Proposition \ref{kveryamp}. Then the last statement follows
by \cite[Theorem 3.3 (b)]{Ful98}.
\end{proof}

As an application, we prove the theorems mentioned earlier.

\begin{proof}[Proof of Theorem \ref{thm4}]
    It follows by combining Theorem \ref{thm:main41}, (\ref{eq:ch1cd}), Propositions \ref{prop:lth2} and \ref{prop:lth3}. In the situations (1)-(3), we have $h^0(\cO_S(\beta))\ge d+1$ because in each case, the morphism $\sigma_{d-1}:\cC^{[d]}\to S^{[d]}$ is surjective.
\end{proof}

\begin{proof}[Proof of Theorem \ref{mt:piceq}]
It suffices to show $\rho(M_{\beta,\chi})\le\rho(S)+1$ in the situations (1)-(3) of Theorem \ref{thm4}. Since $M_{\beta,\chi}$ is normal, $\Pic(M_{\beta,\chi})_\bQ\subset\CH^1(M_{\beta,\chi})_\bQ$ by \cite[Ch. \uppercase\expandafter{\romannumeral2}, Corollary 6.14]{Hartshorne}. Combining Theorem \ref{thm4}, \cite[Corollary 6.3]{Fog73}, and our assumption $h^1(\cO_S)=0$, we have
$$\rho(M_{\beta,\chi})\le\dim\Pic(M_{\beta,\chi})_\bQ\le\dim\CH^1(M_{\beta,\chi})_\bQ\leq\dim\CH^1(S^{[d]})_\bQ\le\rho(S)+1,$$
which completes the proof.
\end{proof}



\subsection{Further discussion and some questions}
We believe that one key ingredient 
\begin{equation} \label{id:hil}
    \CH^1(\cC^{[d]})_\bQ=\CH^1(S^{[d]})_\bQ\oplus \bQ
\end{equation} in our approach should hold more generally. We will discuss possible approaches to dealing with (\ref{id:hil}) and related questions for future directions. 
\subsubsection{A possible approach via relative Hilbert--Chow}
Let $\cC^{(d)}$ be the relative $d$-th symmetric product of the universal curve $\pi:\cC \rightarrow |\beta|$, i.e., $\cC^{(d)}:=\cC_{\pi}^d/\mathfrak{S}_d$ is the quotient of the $d$-th fiber product $$\cC_{\pi}^d:=\cC \times_{|\beta|} \times \cdots \times_{|\beta|} \cC$$
by the permutation action of the symmetric group 
$\mathfrak{S}_d$. 
Then there is a relative Hilbert-Chow morphism  ({\it cf}. \cite[\S 1.2]{Ran16})
\begin{equation*}
    h:\ \cC^{[d]} \rightarrow \cC^{(d)}.
\end{equation*}
Note that $\cC^{(d)}$ is irreducible by Proposition \ref{prop:relhilirr} and the surjectivity of $h$. 
 Denote by $B\subset |\beta|$ the locus of integral curves with at worst nodal singularities. Let $\pi_B:\cC_B \rightarrow B$, $\cC^{[d]}_{B}$, and $\cC_B^{(d)}$ be the base change of $\pi$, $\cC^{[d]}$, and $\cC^{(d)}$, respectively. Since by Proposition \ref{eqdim:cC}, $\pi^{[d]}:\cC^{[d]} \rightarrow |\beta|$ is equidimensional of relative dimension $d$, we have $\codim(\cC^{[d]}\setminus\cC^{[d]}_{B},\cC^{[d]})\geq2$ provided  $\codim(|\beta|\setminus B,|\beta|)\ge 2$, which is true in many cases. Assume this is true from now on and thus \[ \dim \CH^1(\cC^{[d]})_\bQ=\dim \CH^1(\cC^{[d]}_{B})_\bQ.\]
 Moreover, by \cite{Ran05}, $\cC^{[d]}_{B}$ is smooth and thus $\CH^1(\cC^{[d]}_{B})_\bQ=\Pic(\cC^{[d]}_{B})_\bQ$.  Now we only consider the relative Hilbert-Chow morphism $h$ over $B$. It was shown in \cite[Theorem 2.1]{Ran16} that $h$ is the blow-up along a Weil but not $\bQ$-Cartier divisor $D^{[d]}\subset \cC^{(d)}_{B}$, then \[ \CH^1( \cC^{(d)}_{B})_\bQ= \CH^1(\cC^{[d]}_{B})_\bQ, \ \ \ \rho(\cC^{[d]}_{B})=\rho(\cC^{(d)}_{B})
 +1.\]

Parallel to the absolute case, we expect 
that similar results still hold for the relative case. Thus, we ask the following 
\begin{question}\label{Ques:Relativepicard}
    Is it true that $\dim \Pic(\cC^{(d)}_{B})_\bQ=\dim \Pic(\cC_{B})_\bQ$?
\end{question}
\begin{rem}  
   Heuristically, if we view $\cC_{B}$ as a smooth projective curve $C_K$ over the function field $K:=\bC(B)$, then $\cC^{(d)}_{B}$ should be viewed as $d$-th symmetric product $C_K^{(d)}$. As an analogous result over $\bC$ ({\it cf}. \cite[Proposition 2.3]{Mustopa08}),  ``$\NS(C_K^{(d)})=\NS(J_{C_K})\oplus \bZ$'' is expected. Then it is also reasonable to expect $\NS(\cC^{(d)}_{B})=\NS(\cC_{B})\oplus \bZ$. As $\cC^{(d)}_{B}$ has only one Weil divisor which is not Cartier by Ran's work \cite{Ran16}, this suggests a positive answer to Question \ref{Ques:Relativepicard}.
\end{rem}
Now assume the above Question \ref{Ques:Relativepicard} has a positive answer, $h^1(\cO_S)=0$, and $\CH^1(\cC^{[d]})_\bQ$ is finite dimensional. Then by the above discussion, 
\begin{equation*}
    \begin{split}
        \dim \CH^1(\cC^{[d]})_\bQ& =\dim \CH^1(\cC^{[d]}_{B})_\bQ= \dim \Pic(\cC^{[d]}_{B})_\bQ \\ &=\dim \Pic(\cC^{(d)}_{B})_\bQ
 +1=\dim \Pic(\cC_{B})_\bQ+1\\ &=\rho(S)+2.
    \end{split}
\end{equation*}
 So Question \ref{ques:pic} has positive answer in more general cases combined with the lower bound in Theorem \ref{picardnum}. Towards finite dimensionality, we have the following proposition.

\begin{pro}
\label{prop:picfrk}
    Assume $h^1(\cO_S)=0$ and $h^1(\cO_S(\beta+K_S))=0$ (e.g., $\beta$ is nef and big). Then $h^1(\cO_{\cC_{\pi}^n})=0$ and $h^1(\cO_{\cC^{(n)}})=0$ for all $n\in\bZ_{>0}$.
\end{pro}
\begin{proof}
    It is sufficient to show the first assertion since $\cC^{(n)}=\cC_{\pi}^n/\mathfrak{S}_n$. We will show this by induction on $n\in \bZ_{>0}$. The $n=1$ case is easy since $\cC$ is a projective bundle over $S$ by the base-point-freeness of $\beta$ and Proposition \ref{kveryamp}. Assuming $h^1(\cO_{\cC_{\pi}^n})=0$, we need to show $h^1(\cO_{\cC_{\pi}^{n+1}})=0$.
Consider the following diagram
\begin{equation*}
    \begin{tikzcd}
        \cC_{\pi}^{n+1} \arrow[r,hook,"\iota"]\arrow[d,"p"]& \cC_{\pi}^{n} \times S    \arrow[dl,"pr_1"] \\
        \cC_{\pi}^{n}. & 
    \end{tikzcd}
\end{equation*}
The Cartier divisor $\cC_{\pi}^{n+1}\subset\cC_{\pi}^{n} \times S$ corresponds to the line bundle $A_n\boxtimes \cO_S(\beta)$, where $A_n$ is the pullback of $\cO_{|\beta|}(1)$ along $\cC_{\pi}^{n}\to|\beta|$. By the induction hypothesis and the Künneth formula, $h^1(\cO_{\cC_{\pi}^{n} \times S})=0$. By the long exact sequence associated to \[  0\rightarrow (A_n\boxtimes \cO_S(\beta) )^\vee \rightarrow \cO_{\cC_{\pi}^{n} \times S} \rightarrow  \iota_\ast\cO_{ \cC_{\pi}^{n+1}}  \rightarrow 0,\]
we only need to show $H^2(\cC_{\pi}^{n} \times S,(A_n^\vee\boxtimes \cO_S(-\beta)) )=0$. Since $A_n$ and $\cO_S(\beta)$ are nontrivial line bundles with nonzero sections, $h^0(\cO_S(-\beta))=h^0(A_n^\vee)=0$. By Serre duality and our assumption, $h^1(\cO_S(-\beta))=h^1(\cO_S(\beta+K_S))=0$.
Hence the result follows from the Künneth formula.
\end{proof}

By the exponential sequence, Proposition \ref{prop:picfrk} implies that $\Pic( \cC_{\pi}^{d})$ and $\Pic( \cC^{(d)})$ have finite ranks. In particular, their spaces of Cartier divisors with $\bQ$-coefficients are finite dimensional. Combining the resolution of singularities and localization sequence of Chow groups, one can easily see that $\CH^1(\cC^{[d]})_\bQ$ is finite dimensional under the assumption  $h^1(\cO_S)=0$.

\subsubsection{Brill--Noether type loci}
We discuss the geometry of $W^d_{\le k}(\beta)$ $(d \in\bZ_{>0})$ from the perspective of Brill--Noether loci on surfaces, which is very active recently. See the survey \cite{CHN24} for positive rank cases.
 Assume $h^1(\cO_S(\beta))=0$. By the long exact sequence induced by (\ref{ses}), \[ W^d_{\le k}(\beta)=\{ Z\in S^{[d]} |\ h^0(\cI_Z\otimes \cO_S(\beta)) \ge N+1-k\} \]
or equivalently 
\[ W^d_{\le k}(\beta)=\{ Z\in S^{[d]} |\ h^1(\cI_Z\otimes \cO_S(\beta))=\mathrm{ext}^1(\cO_S(-\beta),\cI_Z) \ge d-k \}. \]
We call this locus the $k$-th $\beta$-twisted Brill--Noether locus in $S^{[d]}$, as similar loci have been studied by Bayer--Chen--Jiang \cite{BCJ24}. In \cite{BCJ24} it was shown that the Brill--Noether locus in $S^{[d]}\times S$ is of expected codimension by explicit coordinate calculation. In our situation, for a given affine open subset $\operatorname{Spec}(R)$ of $S^{[d]}$, 
the ideal of $ W^d_{\le k}(\beta)$ on $\operatorname{Spec}(R)$ is generated by $\binom{N+1}{k+1}\cdot\binom{d}{k+1}$ elements in $R$. 
\begin{question}
    Is there a way to prove (\ref{codimcond}) by finding ($d+2$) independent generators among those $\binom{N+1}{k+1}\cdot\binom{d}{k+1}$ elements when $\beta$ is sufficiently positive? 
\end{question}
\subsubsection{Deeper analysis of wall-crossing}
Note that the only essential wall-crossing we used is to connect the moduli spaces $P_{\beta,\chi}(0^+)$ and $P_{\beta,\chi}(\infty)$. It is more ambitious to describe the geometry of wall-crossing \[ P_{\beta,\chi}(\delta_i) \stackrel{p_-}\longrightarrow  P_{\beta,\chi}(w_{i}) \stackrel{p_+}\longleftarrow   P_{\beta,\chi}(\delta_{i-1}) \] at each wall $w_i\in \bQ_{>0}$, where $\delta_i\in(w_{i},w_{i+1})\cap\bQ$ and $\delta_{i-1}\in(w_{i-1},w_i)\cap\bQ$. More precisely, we would like to ask the following 
\begin{question}\label{Ques:WC}
    What are the fibers of $p_\pm$? How to describe the center $Z_{\pm}\subset  P_{\beta,\chi}(w_{i})$? 
\end{question}
Here the center $Z_{\pm}\subset  P_{\beta,\chi}(w_{i})$ is the locus where $p_\pm$ is not isomorphic. A good understanding of Question \ref{Ques:WC} should be helpful in relating the stabilization of relative Hilbert schemes and that of moduli spaces of one-dimensional sheaves.

\section{Reducible moduli spaces on certain surfaces of general type}
\label{sec:other}


In this section, we give examples where the asymptotic irreducibility fails for certain surfaces, which indicates that more complicated phenomena may occur for surfaces of general type.
We make the following assumption.
\begin{assu}
\label{assu:irrd}
The following conditions hold:
    \begin{enumerate}
    \item $h^1(\cO_S) \ge 1$.
    \item A general curve in $|K_S|\not=\emptyset$ is smooth and connected.
    \item \label{cond:n} There exists an integer $n\ge 2$ such that $h^1(\cO_S(nK_S))=0$ and that a general curve in $|nK_S|$ is smooth, connected.
\end{enumerate}
\end{assu}

Our proof of reducibility is the ``opposite'' of that of \cite[Theorem 2.3]{MS23}: instead of proving that $\dim h_{\beta,\chi}^{-1}(nC)$ is not too big, we show it is not too small.
\begin{thm}\label{counter} Under Assumption \ref{assu:irrd},
    the moduli space $M_{\beta_n,1}$ is reducible for $\beta_n=nK_S$ with $n$ as in the condition (\ref{cond:n}).
\end{thm}
\begin{proof}
By the Riemann--Roch formula and Serre duality,
\begin{equation}
\label{dimks}
    \begin{aligned}
        \dim|K_S|&=\chi(K_S)+h^1(K_S)-h^2(K_S)-1\\
        &=\left(\frac{K_S(K_S-K_S)}{2}+\chi(\cO_S)\right)+h^1(\cO_S)-h^0(\cO_S)-1\\
        &=\chi(\cO_S)+h^1(\cO_S)-2.
    \end{aligned}
\end{equation}
For a general $C\in|K_S|$, which is smooth and connected by the condition (2) of Assumption \ref{assu:irrd}, the fiber $h_{n}^{-1}(nC)$ of the Hilbert--Chow morphism $h_n:=h_{\beta_n,1}$ contains the moduli space of stable vector bundles of rank $n$ on $C$. Hence by the deformation-obstruction theory of moduli of sheaves on curves ({\it cf}. \cite[Corollary 4.5.5]{HL10}) and the adjunction formula,
\begin{equation}
\label{dimfib}
    \begin{aligned}
     \dim h_{n}^{-1}(nC)&\geq n^2(g(C)-1)+1\\
     &=\frac{n^2}{2}(K_S(K_S+K_S))+1=\beta_n^2+1.
    \end{aligned}
\end{equation}
We denote by $U_n$ the open subset of $|\beta_n|$ consisting of smooth curves, and by $W_n$ the complement of $U_n$ in $|\beta_n|$.
By (\ref{dimks}) and (\ref{dimfib}), we have
\begin{equation}
\label{ineq:irrd}
\begin{aligned}
    \dim h_n^{-1}(W_n)&\geq \dim|K_S|+\dim h_n^{-1}(nC)\\
    &\geq(\chi(\cO_S)+h^1(\cO_S)-2)+(\beta_n^2+1)\\
    &=\dim h_n^{-1}(U_n)+(h^1(\cO_S)-1),
\end{aligned}
\end{equation}
where the last identity follows from (\ref{dimsmsupp}) combined with the condition (\ref{cond:n}) of Assumption \ref{assu:irrd}.
Now we finish the proof by contradiction. 
Suppose the assertion is false, i.e.,  $M_{\beta_n,1}$ is irreducible. Then
$$\dim h_n^{-1}(W_n)<\dim M_{\beta_n,1}=\dim h_n^{-1}(U_n)$$
since $h_n^{-1}(U_n)\not=\emptyset$ is open and $h_n^{-1}(W_n)$ is a proper closed subset of $M_{\beta_n,1}$,
which contradicts (\ref{ineq:irrd}).
\end{proof}
We give some explicit examples satisfying Assumption \ref{assu:irrd}.

\begin{exmp}
Let $X$ be an abelian threefold with a very ample divisor $L$. Then a general $S\in|L|$ is a surface whose canonical divisor is very ample by the adjunction formula. Hence $S$ satisfies the conditions (2) and (3) in Assumption \ref{assu:irrd} by Bertini's theorem and Kodaira vanishing. The condition (1) follows from the Lefschetz hyperplane theorem.
\end{exmp}

\begin{exmp}
Let $A$ be a principally polarized abelian surface and let $\Theta$ be the theta divisor on $A$. Consider the double cover $S \rightarrow A$ branched along a smooth curve $C_m\in |2m\Theta|$ for $m\in \bZ_{>0}$. Then $K_S$ is ample and $h^1(\cO_S)=h^1(\cO_A)\geq1$. When $m$ and $n$ are large, Assumption \ref{assu:irrd} is satisfied.
\end{exmp}
Note that the argument in Proposition \ref{counter} fails if the curve $C$ lies in $|nK_S|$ for $n\ge 2$. From the general philosophy of asymptotic phenomena in the higher rank case, we still ask if the following question is true:
\begin{question}
  Is $M_{\beta,1}$ irreducible when Assumption \ref{assu:irrd} fails and $\beta$ is sufficiently positive?   
\end{question}

\bibliographystyle{alpha}
\bibliography{main}
\end{document}